# Generalized Fractional Order Derivatives, Its Properties and Applications

*Ali KARCI*
İnönü University
Faculty of Engineering
Department of Computer Engineering
44280 Malatya / Turkey

## Abstract

The concept of fractional order derivative can be found in extensive range of many different subject areas. For this reason, the concept of fractional order derivative should be examined. After giving different methods mostly used in engineering and scientific applications, the omissions or errors of these methods will be discussed in this study. The mostly used methods are Euler, Riemann-Liouville and Caputo which are fractional order derivatives. The applications of these methods to constant and identity functions will be given in this study. Obtained results demonstrated that all of three methods have errors and deficiencies. In fact, the obtained results demonstrated that the methods given as fractional order derivatives are curve fitting methods or curve approximation methods. In this paper, we redefined fractional order derivative (FOD) by using classical derivative definition and L'Hospital method, since classical derivative definition concluded in indefinite limit such as 0/0. The obtained definition is same as classical derivative definition in case of fractional order is equal to 1. This implies that definitions and theorems in this paper sound and complete.

After defining the FOD concept, the applications of FOD for polynomial, exponential, trigonometric and logarithmic functions were handled in this study. The properties of FOD for positive monotonic increasing/decreasing functions were demonstrated in this paper. Another important point is that the relation between derivatives and complex functions were also verified in this paper.

**Keywords.** Fractional Order Derivatives, Euler Method, Riemann-Liouville Method, Caputo Method, Complex Functions.

## 1 Introduction

The FODs are different approaches instead of classical derivatives. There are a lot of studies on this subject. The most of these studies have been used Euler, Riemann-Liouville and Caputo FODs. Due to this case, this study focused on Euler, Riemann-Liouville and Caputo FODs.

Some of studies on the FODs can be summarized as follows.

A minimization problem with a Lagrangian that depends on the left Riemann–Liouville FOD was considered in [1] such as finite differences, as a subclass of direct methods in the calculus of variations, consist in discretizing the objective functional using appropriate approximations for derivatives that appear in the problem. There is a study on fractional extensions of the classical Jacobi polynomials [2], and fractional order Rodrigues' type representation formula. By means of the Riemann–Liouville operator of fractional calculus, new g-Jacobi functions were defined, some of their properties were given and compared with the corresponding properties of the classical Jacobi polynomials [2]. There is another study on the discussion of theory of fractional powers of operators on an arbitrary Frechet space, and the authors of this study obtained multivariable fractional integrals and derivatives defined on certain space of test functions and generalized functions [3].

Differential equations of fractional order appear in many applications in physics, chemistry and engineering [4]. There is a requirement for an effective and easy-to-use method for solving such equations. Bataineh et al used series solutions of the fractional differential equations using the homotopy analysis method [4]. Many recently developed models in areas like viscoelasticity, electrochemistry, diffusion processes, etc. are formulated in terms of derivatives (and integrals) of fractional (non-integer) order [5]. There is a collection of numerical algorithms for the solution of the various problems arising in derivatives of fractional order [5].

The fractional calculus is used to model various different phenomena in nature but due to the non-local property of the fractional derivatives, it still remains a lot of improvements in the present numerical approaches [6]. There are some approaches based on piecewise interpolation for fractional calculus, and some improvement based on the Simpson method for the fractional differential equations [6]. There is a study on FOD learning control including many theoretical and experimental results, and these results shown the improvement of transient and steady-state performances [7].

Recently, many models are formulated in terms of FODs, such as in control processing, viscoelasticity, signal processing, and anomalous diffusion, and authors of this paper studied the important properties of the Riemann-Liouville derivative, one of mostly used FODs [8]. The philosophy of integer order sliding mode control is valid also for the systems represented by fractional order operators [9].

Here is just a small part of the work in this area is given. Except these, there are numerous studies. This area of work can be divided into two groups such as mathematical theory of FODs and the applications of FODs. There are important errors and deficiencies in both groups of studies. The sources and reasons for these errors and deficiencies can be summarized as follows:

a) It was assumed that the order of derivative is integer up to a specific step. While derivation reached to that step, the order of derivative was converted in to real number. This process concluded in involvement of gamma functions due to the assumption of order of derivative as integer. This is a deficiencies and error. For example, $f(x)=x^k$ and $n<k$. $n,k\in Z^+$, $n^{th}$ order derivative of f(x) is $f^{(n)}(x)$.

$$f^{(n)}(x)=k(k-1)(k-2)\ldots(k-n+1)x^{k-n}$$

and $(n+1)^{th}$ is assumed as FOD. If there is not any such assumption, the process was not concluded in involvement of gamma functions. There is a generalization that FOD and classical derivative have been assumed as same up to a specific step. This is the source of deficiencies and errors. Assume that the order of derivative is $\frac{\beta}{\delta}$. The coefficient of $x^{k-n}$ variable is

$$k\left(k-\frac{\beta}{\delta}\right)\left(k-2\frac{\beta}{\delta}\right)\ldots\left(k-(n-1)\frac{\beta}{\delta}\right).$$

This is not a gamma function.

b) Another important point is that powers of variables decreased as integer during derivation process. There was another assumption in FOD regarded the powers of variables. The coefficients of variables were determined with gamma functions, and powers declining of variables were also regarded as integer up to a specific step, and after reaching that step, power decreasing was assumed as real number. For example,

$f(x)=x^k$ and $n<k$. $n,k\in Z^+$, $n^{th}$ order derivative of f(x) is $f^{(n)}(x)$, then

$$f^{(n)}(x)=k\left(k-\frac{\beta}{\delta}\right)...\left(k-(n-1)\frac{\beta}{\delta}\right)x^{k-n\frac{\beta}{\delta}}$$

The obtained equation for $f^{(n)}(x)$ is not a derivative. Since it is straightforward process of decreasing power of variable step-by-step and it is also product of powers of variable at each step as coefficient.

In this study, it is focused on the shortcomings and wrong points involved in the methods of Euler, Riemann-Liouville and Caputo for FODs. Especially, the FODs of constant and identity functions will be obtained in this study. Euler and Riemann-Liouville methods were yielded shortcomings and errors in results for constant functions, and on contrary, Caputo method was yielded in correct result for constant function. The methods have not provided accurate results for identity function. So, all three methods are curve fitting or curve approximation methods.

## 2 Deficiencies of Old Fractional Order Derivative Methods

There are different methods and approximations for FODs since 1730. There is a common property related to these methods and approximations. The order of derivative as integer caused gamma function involvement in most of the methods and approximations. We used three most popular methods in this study.

a) L.Euler (1730) method:

$$\frac{d^n x^m}{dx^n}=m(m-1)...(m-n+1)x^{m-n} \tag{1}$$

By using gamma function, this method can be rewritten as in Eq. 2.

$$\frac{d^n x^m}{dx^n}=\frac{\Gamma(m+1)}{\Gamma(m-n+1)}x^{m-n} \tag{2}$$

The definition of gamma function is seen in Eq.3.

$$\Gamma(z)=\int_0^\infty e^{-t}t^{z-1}dt \tag{3}$$

b) Riemann-Liouville method: The method of Riemann-Liouville is given in Eq.4.

$${}_aD_t^\alpha f(t)=\frac{1}{\Gamma(n-\alpha)}\left(\frac{d}{dt}\right)^n\int_a^t\frac{f(v)dv}{(t-v)^{\alpha-n+1}} \tag{4}$$

c) Caputo (1967) method: The method of Caputo is seen in Eq.5.

$${}_a^C D_t^\alpha f(t)=\frac{1}{\Gamma(\alpha-n)}\int_a^t\frac{f^{(n)}(v)dv}{(t-v)^{\alpha+1-n}} \tag{5}$$

### 2.1 Fractional Order Derivative of $f(x)=cx^0$

These three methods for FODs are mostly used methods. Due to this case, we will investigate the fractional order derivative of each method in this section.

First of all, the results of Euler, Riemann-Liouville and Caputo methods for $f(x)=cx^0$ are illustrated in Eq.6, Eq.8, and Eq.9. c is a constant, n=1 and $\alpha=\frac{2}{3}$ for all methods.

*Euler:*

$$\frac{d^n x^m}{dx^n}=\frac{\Gamma(m+1)}{\Gamma(m-n+1)}x^{m-n}$$
$$=\frac{d^{\frac{2}{3}}x^0}{dx^{\frac{2}{3}}} \tag{6}$$
$$=\frac{\Gamma(1)}{\Gamma\left(\frac{1}{3}\right)}cx^{-\frac{2}{3}}$$

$\Gamma(1)=1$ and $\Gamma\left(\frac{1}{3}\right)$ is a transcendental number such as π and base of natural logarithm e= 2,7134884828. Since

$$\Gamma\left(\frac{1}{3}\right)=\int_0^\infty t^{-\frac{2}{3}}e^{-t}dt$$
$$=3\sum_{i=0}^{\infty}\frac{3^i e^{-t} t^{\left(\frac{3i+1}{3}\right)}}{\prod_{j=0}^{i}(3j+1)} \tag{7}$$

and Eq.7 illustrates that it cannot be computed in a simple way. $f(x)=cx^0$ is a constant function and its change with respect to the change in x is zero. However, Eq.6 depicts that the result obtained from Euler method is not zero and it depends on variable x. The two conditions must be taken in care.

i) If $\Gamma\left(\frac{1}{3}\right)\neq 0$ and is finite, then the FOD of $f(x)=cx^0$ is concluded in dependence on variable x. The change in a constant is zero in all cases.

ii) If $\Gamma\left(\frac{1}{3}\right)=0$, then $\frac{d^n x^m}{dx^n}=\frac{\Gamma(1)}{\Gamma\left(\frac{1}{3}\right)}x^{-\frac{2}{3}}$ for n=1 and $\alpha=\frac{2}{3}$ is indefinite. The obtained results in both cases are inconsistent.

*Riemann-Liouville method*:

Assume that n=1 and $\alpha=\frac{2}{3}$

$$
\begin{aligned}
{}_aD_t^{\alpha} f(t) &= \frac{1}{\Gamma(n-\alpha)}\left(\frac{d}{dt}\right)^n \int_a^t \frac{f(v)dv}{(t-v)^{\alpha-n+1}} \\
&= {}_aD_t^{\frac{2}{3}} f(t) \\
&= \frac{1}{\Gamma\left(\frac{1}{3}\right)} \frac{d}{dt} \int_a^t \frac{cdv}{(t-v)^{\frac{2}{3}}} \\
&= \frac{c}{\Gamma\left(\frac{1}{3}\right)}\left(-\frac{1}{(t-a)^{\frac{2}{3}}}\right) \\
&\neq 0
\end{aligned} \tag{8}
$$

The obtained result is inconsistent, since the result is a function of x. However, initial function is a constant function.

*Caputo method*:

Assume that n=1 and $\alpha = \frac{2}{3}$

$$
\begin{aligned}
{}_a^CD_t^{\alpha} f(t) &= \frac{1}{\Gamma(\alpha-n)} \int_a^t \frac{f^{(n)}(v)dv}{(t-v)^{\alpha+1-n}} \\
&= \frac{1}{\Gamma\left(-\frac{1}{3}\right)} \int_a^t \frac{0dv}{(t-v)^{\frac{2}{3}+1-1}} \\
&= 0
\end{aligned} \tag{9}
$$

The result of Caputo method is consistent.

The Euler and Riemann-Liouville methods do not work for constant functions as seen in Eq.6 and Eq.8. There is no change in constant function. If there is any change in constant function, it is not a constant function. On contrary, any order derivative of constant function is zero with respect to Caputo method.

### 2.2 Fractional Order Derivative of f(x)=x

Fractional order derivatives of identity function obtained with respect to Euler, Riemann-Liouville, Caputo methods in this section. Assume that n=1 and $\alpha = \frac{2}{3}$.

*Euler method*:

$$
\begin{aligned}
\frac{d^{\frac{2}{3}} x^1}{dx^{\frac{2}{3}}} &= \frac{\Gamma(1+1)}{\Gamma\left(1-\frac{2}{3}+1\right)} x^{1-\frac{2}{3}} \\
&= \frac{\Gamma(2)}{\Gamma\left(\frac{4}{3}\right)} x^{\frac{1}{3}} \\
&\neq 1
\end{aligned} \tag{10}
$$

Classical definition of derivative is

$$\lim_{h \to 0} \frac{f(x+h)-f(x)}{(x+h)-x}.$$

It can be seen from definition, the ratio of the change in dependent variable to the change in independent is always 1 (one). In this case, the derivative must be 1 in any FOD. However, FOD of identity function with respect to Euler method is different from 1. This means that Euler method yielded in an inconsistent result.

*Riemann-Liouville method:*

$$\begin{aligned}
{}_aD_t^{\alpha} f(t) &= \frac{1}{\Gamma(n-\alpha)}\left(\frac{d}{dt}\right)^n \int_a^t \frac{f(v)dv}{(t-v)^{\alpha-n+1}} \\
&= {}_aD_t^{\frac{2}{3}} f(t) \\
&= \frac{1}{\Gamma\left(\frac{1}{3}\right)} \frac{d}{dt} \int_a^t \frac{xdx}{(t-x)^{\frac{2}{3}-1+1}} \\
&= \frac{1}{\Gamma\left(\frac{1}{3}\right)} \frac{d}{dt} \int_a^t \frac{xdx}{(t-x)^{\frac{2}{3}}} \\
&= \frac{1}{\Gamma\left(\frac{1}{3}\right)} \left(3a(t-a)^{\frac{2}{3}} + \frac{9}{4}(t-a)^{\frac{4}{3}}\right) \\
&\neq 1
\end{aligned} \tag{11}$$

The obtained result with respect to Riemann-Liouville method is an inconsistent result.

*Caputo method*:

$$\begin{aligned}
{}_a^C D_t^{\alpha} f(t) &= \frac{1}{\Gamma(\alpha-n)} \int_a^t \frac{f^{(n)}(v)dv}{(t-v)^{\alpha+1-n}} \\
&= \frac{1}{\Gamma\left(-\frac{1}{3}\right)} \int_a^t \frac{dv}{(t-v)^{\frac{2}{3}+1-1}} \\
&= \frac{1}{\Gamma\left(-\frac{1}{3}\right)} \left(3(t-a)^{\frac{1}{3}}\right) \\
&\neq 1
\end{aligned} \tag{12}$$

The obtained result with respect to Caputo method is also an inconsistent result.

The deficiencies in all three methods for identity function (f(x)=x) and in two methods (Euler, Riemann-Liouville) for f(x)=c are due to the attention to coefficients and powers of functions during taking derivatives. Another important point in the process of differentiation is that integer coefficients and powers in obtaining formula for differentiation are assumed. After this process, coefficients and powers are assumed as real numbers in a specific step. The formula obtained in this way is not a FOD, they can be regarded as curve fitting.

Euler, Riemann-Liouville and Caputo methods do not work as FOD for identity function (Eq.10, Eq.11 and Eq.12). The assumptions of coefficients and powers makes formula obtained by Euler, Riemann-Liouville, Caputo methods are curve fitting formula.

### 2.3 Applications of Old FOD Methods

In this study, the most important three methods for FODs were investigated for constant and identity functions. The derivative of any order for $f(x)=cx^0$ is zero, however, Fig.1 and Fig.2 depict that the derivative of $f(x)=cx^0$ is not zero for Euler and Riemann-Liouville methods.

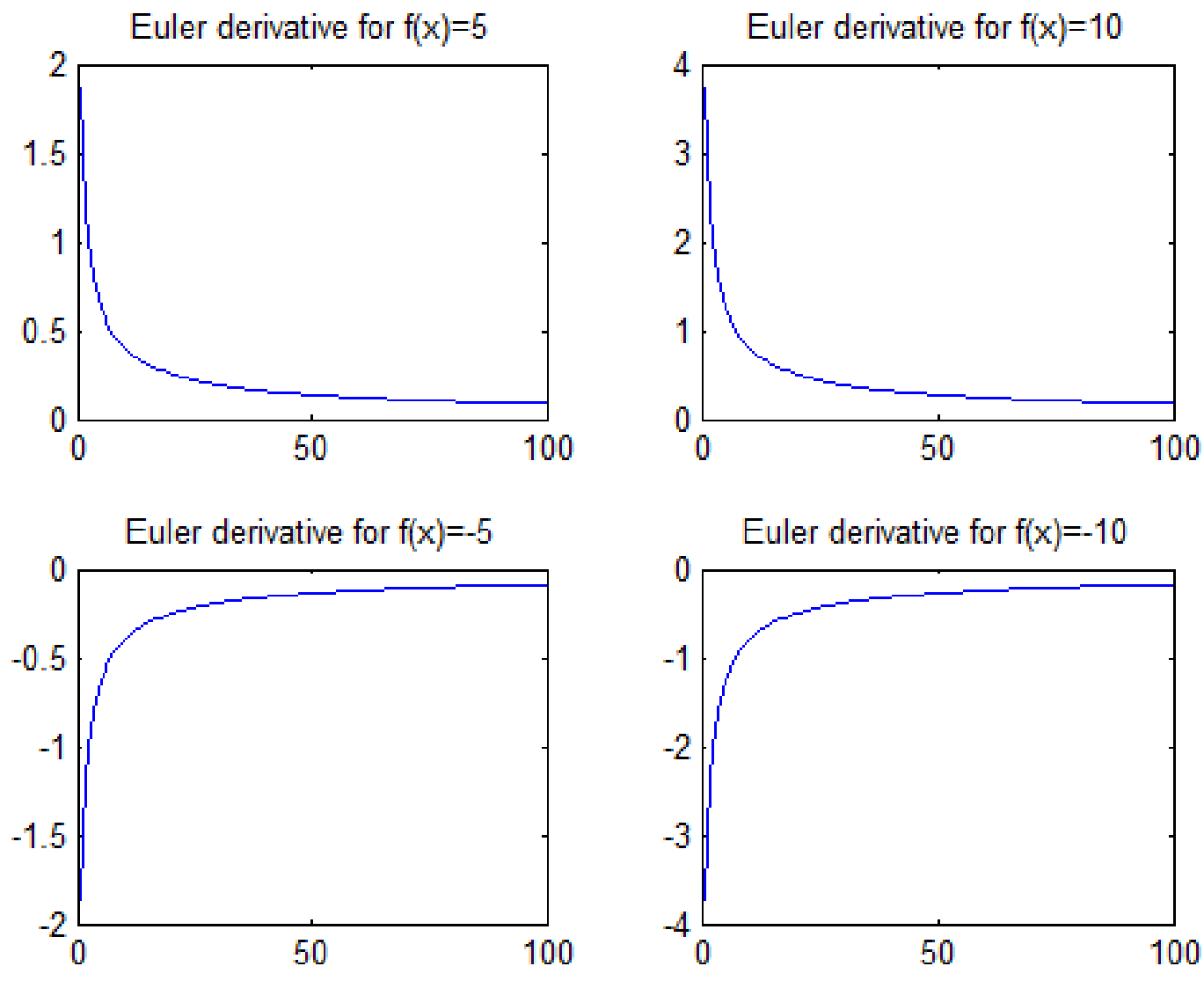

**Figure 1.** Figures of FODs by using Euler method for $f(x)=cx^0$ (c=5, c=10, c=-5, c=-10).

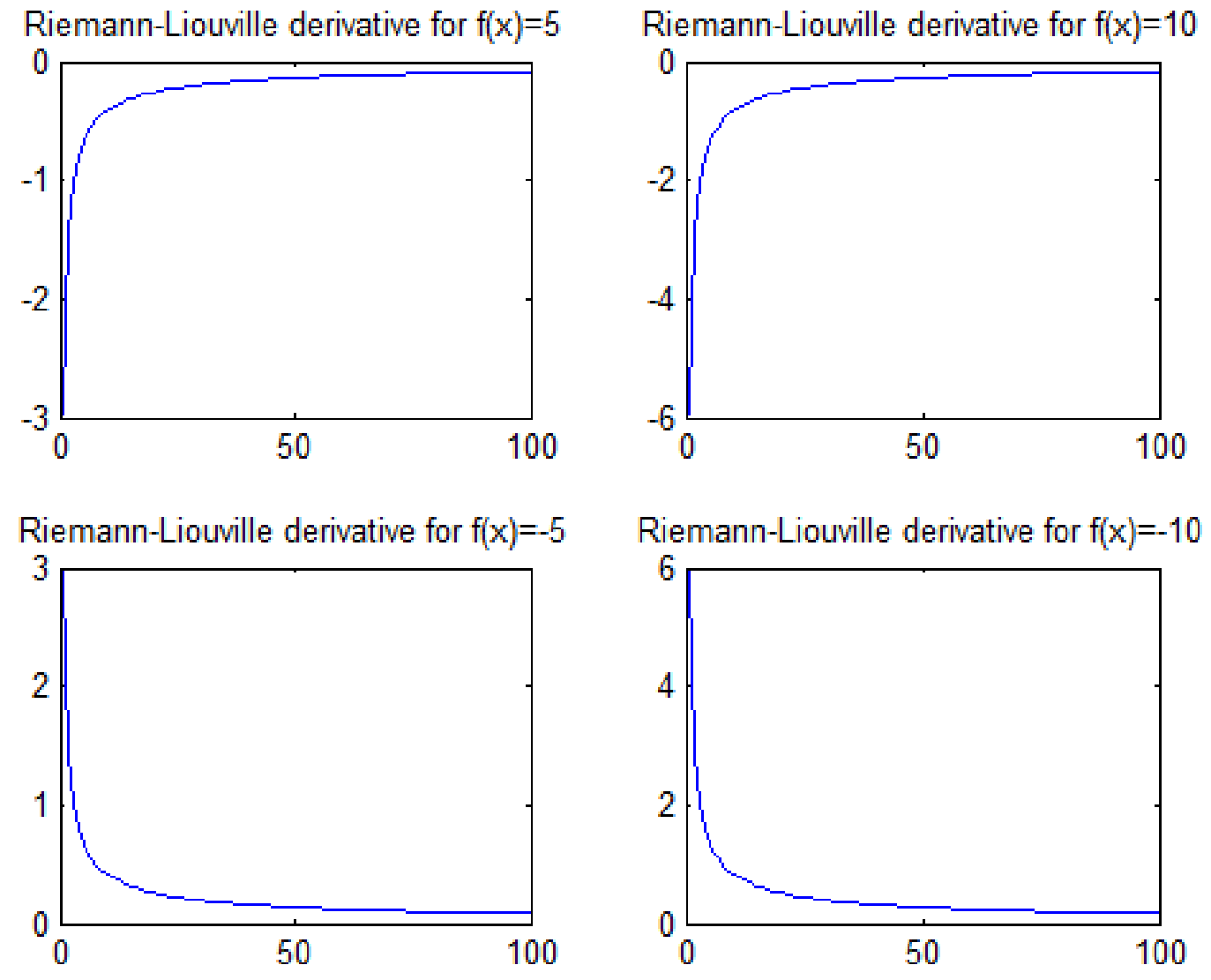

**Figure 2.** Figures of FODs by using Riemann-Liouville method for $f(x)=cx^0$ (c=5, c=10, c=-5, c=-10).

Fig.1 depicts the FODs by using Euler method for f(x)=5, f(x)=10, f(x)=-5 and f(x)=-10. The derivatives must be zero for four cases due the constant functions. However, the obtained results are not zero; this is an evidence for deficiency of Euler method.

Fig.2 depicts the FODs by using Riemann-Liouville method for f(x)=5, f(x)=10, f(x)=-5 and f(x)=-10. The derivatives must be zero for four cases due the constant functions. However, the obtained results are not zero; this is an evidence for deficiency of Riemann-Liouville method.

**Table 1.** Euler and Riemann-Liouville methods Matlab source code.

```
for i=1:100
  Euler(i)=5*(1/gamma(1/3)*(i^(-2/3)));
  Riemann(i)=5*((1/gamma(1/3))*(-1/((x(i)-
0.5)^(2/3))));
end
subplot(2,2,1) plot(x,Euler);
title('Euler derivative for f(x)=5');
 for i=1:100
  Euler(i)=10*(1/gamma(1/3)*(i^(-2/3)));
  Riemann(i)=10*((1/gamma(1/3))*(-1/((x(i)-
0.5)^(2/3))));
end
subplot(2,2,2) plot(x,Euler);
title('Euler derivative for f(x)=10');
 for i=1:100
  Euler(i)=-5*(1/gamma(1/3)*(i^(-2/3)));
  Riemann(i)=-5*((1/gamma(1/3))*(-1/((x(i)-
0.5)^(2/3))));
end
subplot(2,2,3) plot(x,Euler);
title('Euler derivative for f(x)=-5');
for i=1:100
  Euler(i)=-10*(1/gamma(1/3)*(i^(-2/3)));
  Riemann(i)=-10*((1/gamma(1/3))*(-1/((x(i)-
0.5)^(2/3))));
end
subplot(2,2,4) plot(x,Euler); title('Euler derivative for
f(x)=-10');
 % Riemann-Liouvill
subplot(2,2,1) plot(x,Riemann);
title('Riemann-Liouville derivative for f(x)=5');
subplot(2,2,2) plot(x,Riemann);
title('Riemann-Liouville derivative for f(x)=10');
subplot(2,2,3) plot(x,Riemann);
title('Riemann-Liouville derivative for f(x)=-5');
subplot(2,2,4) plot(x,Riemann);
title('Riemann-Liouville derivative for f(x)=-10');
```

**Table 2.** The source code of Euler, Riemann-Liouville and Caputo FODs for identity function.

```
x=1:100;
a=2.5;
for i=1:100
  Euler(i)=(gamma(2)/gamma(4/3))*(i^(1/3));
  Riemann(i)=(1/gamma(1/3))*(3*a*(i-a)^(2/3)+(9/4)*(i-a)^(4/3));
  Caputo(i)=(1/gamma(1/3))*(3*(i-a)^(1/3));
end
subplot(2,2,1) plot(x,Euler);
title('Euler derivative for f(x)=x');
subplot(2,2,2)  plot(x,Riemann);
title('Riemann-Liouville derivative for f(x)=x');
subplot(2,2,3) plot(x,Caputo);
title('Caputo derivative for f(x)=x');
```

These figures depict that both methods do not work for constant functions, since there is no change in constant. It is assumed in Euler method that coefficient and power are integer numbers up to a specific derivative step. This is a source for deficiency.

It is assumed that the integral process is a derivative process with negative order derivative.  This causes deficiency and the same assumptions for Euler method are also accepted in Riemann-Liouville method. Both of these cases cause the deficiency.

Fig.3 depicts the results of three methods (Euler, Riemann-Liouville, Caputo). The rate of change in the dependent variable is equal to the rate of change in the independent variable, so their ratio must be equal to 1. However, Fig. 3 depicts that the results for three methods are not equal to 1. This means that three methods do not work for identity function.

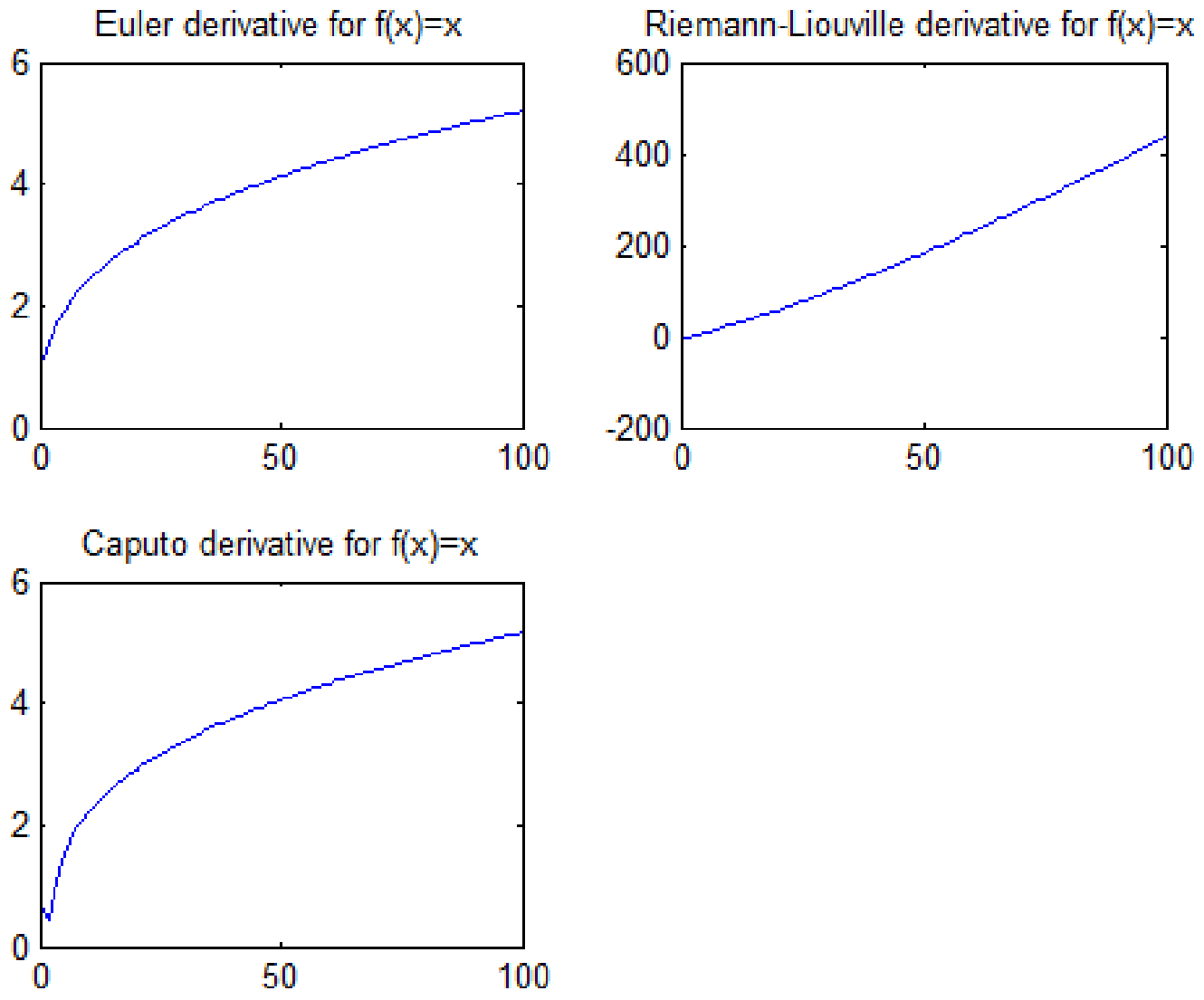

**Figure 3.** The results of Euler, Riemann-Liouville and Caputo FODs.

## 3 New Approach for Fractional Order Derivative

The meaning of derivative is the rate of change in the dependent variable versus the changes in the independent variables. At this aim, the derivative of $f(x)=cx^0$ is that

$$\lim_{h\to 0}\frac{f(x+h)-f(x)}{(x+h)-x}=\lim_{h\to 0}\frac{c-c}{h}=0$$

In the case of identity function is

$$\lim_{h\to 0}\frac{f(x+h)-f(x)}{(x+h)-x}=\lim_{h\to 0}\frac{(x+h)-x}{(x+h)-x}=1$$

So, the definition for fractional order derivative can be considered as follow.

***Definition 1**: f(x):R→R is a function, $\alpha\in R$ and the FOD can be considered as follows.*

$$f^{(\alpha)}(x)=\lim_{h\to 0}\frac{f^{\alpha}(x+h)-f^{\alpha}(x)}{(x+h)^{\alpha}-x^{\alpha}}$$

In the case of very small value of h, the limit in the Definition 1 concluded in indefinite limit.

$$f^{(\alpha)}(x)=\lim_{h\to 0}\frac{f^{\alpha}(x+h)-f^{\alpha}(x)}{(x+h)^{\alpha}-x^{\alpha}}=\frac{0}{0}$$

In this case, the method used for indefinite limit (such as L'Hospital method) can be used, and the FOD can be redefined as follows.

This new approach can be considered for some specific cases.

**Case 1.** f(x)=cx$^0$

$$f^{(\alpha)}(x) = \lim_{h \to 0} \frac{f^{\alpha}(x+h) - f^{\alpha}(x)}{(x+h)^{\alpha} - x^{\alpha}}$$

$$= \lim_{h \to 0} \frac{c - c}{(x+h)^{\alpha} - x^{\alpha}} = 0$$

**Case 2.** f(x)=x

$$f^{(\alpha)}(x) = \lim_{h \to 0} \frac{f^{\alpha}(x+h) - f^{\alpha}(x)}{(x+h)^{\alpha} - x^{\alpha}}$$

$$= \lim_{h \to 0} \frac{(x+h)^{\alpha} - x^{\alpha}}{(x+h)^{\alpha} - x^{\alpha}}$$

$$= 1$$

**Case 3.** f(x) is any function and α>0 (α=2 and f(x)=ax+b)

$$f^{(\alpha)}(x) = \lim_{h \to 0} \frac{f^{\alpha}(x+h) - f^{\alpha}(x)}{(x+h)^{\alpha} - x^{\alpha}}$$

$$= \lim_{h \to 0} \frac{(a(x+h)+b)^2 - (ax+b)^2}{(x+h)^2 - x^2}$$

$$= \frac{a^2 x + ab}{x}$$

**Case 4.** f(x) is any function and α<0 (α=-2 and f(x)=ax+b)

$$f^{(\alpha)}(x) = \lim_{h \to 0} \frac{f^{\alpha}(x+h) - f^{\alpha}(x)}{(x+h)^{\alpha} - x^{\alpha}}$$

$$= \lim_{h \to 0} \frac{(a(x+h)+b)^{-2} - (ax+b)^{-2}}{(x+h)^{-2} - x^{-2}}$$

$$= a\left(\frac{x}{ax+b}\right)^3$$

The definition 1 is a classical definition of derivative, and it has indefinite limit such as $\frac{0}{0}$, while h=0. In this case, Definition 1 can be rephrased as seen in Definition 2.

***Definition 2**: Assume that f(x):R→R is a function, α∊R and L(.)be a L'Hospital process. The FOD of f(x) is*

$$f^{(\alpha)}(x) = \lim_{h \to 0} L\left(\frac{f^{\alpha}(x+h) - f^{\alpha}(x)}{(x+h)^{\alpha} - x^{\alpha}}\right)$$

$$= \lim_{h \to 0} \frac{\dfrac{d\left(f^{\alpha}(x+h) - f^{\alpha}(x)\right)}{dh}}{\dfrac{d((x+h)^{\alpha} - x^{\alpha})}{dh}}$$

***Theorem 1****:Assume that  f(x):R→R is a continue function and $\alpha \in R$. The limit*

$$M = \lim_{h \to 0} L\left(\frac{f^{\alpha}(x+h) - f^{\alpha}(x)}{(x+h)^{\alpha} - x^{\alpha}}\right)$$

*exists and is finite.*

**Proof**: $f^{(\alpha)}(x+h)-f^{(\alpha)}(x)<\infty$ and $(x+h)^{(\alpha)}-x^{(\alpha)}<\infty$.  $f^{(\alpha)}(x+h)-f^{(\alpha)}(x)\neq 0$ and $(x+h)^{(\alpha)}-x^{(\alpha)}\neq 0$. Then

$$\lim_{h \to 0} L\left(\frac{f^{\alpha}(x+h) - f^{\alpha}(x)}{(x+h)^{\alpha} - x^{\alpha}}\right) < \infty$$

and

$$\lim_{h \to 0} L\left(\frac{f^{\alpha}(x+h) - f^{\alpha}(x)}{(x+h)^{\alpha} - x^{\alpha}}\right) = M < \infty \blacksquare$$

The definition of derivative has been handled for α=1 uptil today. While α=1 for definition 2, the obtained results are same to results of classical derivative.

***Theorem 2****: Assume that f(x):R→R is a function, $\alpha \in R$ and α=1, then*

$$\frac{df(x)}{dx} = \lim_{h \to 0} L\left(\frac{f^{\alpha}(x+h) - f^{\alpha}(x)}{(x+h)^{\alpha} - x^{\alpha}}\right).$$

**Proof**: The fractional order derivative of f(x):R→R for α=1 is

$$\frac{df(x)}{dx} = \lim_{h \to 0} \frac{f(x+h) - f(x)}{(x+h) - x}$$

and this derivative implies that

$$\lim_{h \to 0} L\left(\frac{f(x+h) - f(x)}{(x+h) - x}\right) \blacksquare$$

The FOD definition can be demonstrated that it obtained same results as classical derivative definition for α=1. In the subsequent sections, the equality of both methods was illustrated for polinomial, exponential, trigonometric and logarithmic functions.

### 3.1 Equality for Polinomial Functions

For example, $f(x)=x^n$ and α=1, the derivative of function f(x) is

$$\frac{df(x)}{dx} = \frac{(x+h)^n - x^n}{(x+h)-x}$$

$$= \frac{x^n + \binom{n}{1}x^{n-1}h + \binom{n}{2}x^{n-2}h^2 + .... + h^n - x^n}{h}$$

$$= \frac{\binom{n}{1}x^{n-1}h}{h}$$

$$= nx^{n-1}$$

and the FOD is

$$f^{\alpha} = \lim_{h \to 0} L\left( \frac{f^{\alpha}(x+h) - f^{\alpha}(x)}{(x+h)^{\alpha} - x^{\alpha}} \right)$$

$$= \lim_{h \to 0} L\left( \frac{(x+h)^n - x^n}{(x+h)-x} \right)$$

$$= \lim_{h \to 0} L\left( \frac{x^n + \binom{n}{1}x^{n-1}h + \binom{n}{2}x^{n-2}h + ... + h^n}{h} \right)$$

$$= \lim_{h \to 0} \frac{\binom{n}{1}x^{n-1} + \binom{n}{2}x^{n-2}2h + ... + nh^{n-1}}{1}$$

$$= nx^{n-1}$$

### 3.2 Equality for Exponential Functions

For example the exponential function f(x)=$e^x$ and α=1. The derivative of f(x) is

$$\frac{df(x)}{dx} = e^x .$$

The FOD is

$$f^{\alpha}(x) = \lim_{h \to 0} L\left( \frac{e^{x+h} - e^x}{(x+h)-x} \right)$$

$$= \lim_{h \to 0} \frac{e^{x+h}}{1} = e^x$$

### 3.3 Equality for Trigonometric Functions

The derivatives of sinx, cosx, tanx and cotx for α=1 are

$$\frac{df(x)}{dx}=\begin{cases}\cos x & f(x)=\sin x\\ -\sin x & f(x)=\cos x\\ 1+\tan^2 x & f(x)=\tan x\\ -(1+\cot^2 x) & f(x)=\cot x\end{cases}$$

The FODs of trigonometric functions sinx, cosx, tanx and cotx for $\alpha$=1 are

$$\begin{aligned}\frac{d^\alpha f(x)}{dx^\alpha} &= \lim_{h\to 0} L\left(\frac{\sin(x+h)-\sin x}{(x+h)-x}\right)\\ &= \lim_{h\to 0} L\left(\frac{\sin x\cosh+\sinh\cos x-\sin x}{h}\right)\\ &= \lim_{h\to 0}\left(\frac{-\sin x\sinh+\cosh\cos x}{1}\right)\\ &= \cos x\end{aligned}$$

$$\begin{aligned}\frac{d^\alpha f(x)}{dx^\alpha} &= \lim_{h\to 0} L\left(\frac{\cos(x+h)-\cos x}{(x+h)-x}\right)\\ &= \lim_{h\to 0} L\left(\frac{\cos x\cosh-\sinh\sin x-\sin x}{h}\right)\\ &= \lim_{h\to 0}\left(\frac{-\cos x\sinh-\cosh\sin x}{1}\right)\\ &= -\sin x\end{aligned}$$

$$\begin{aligned}\frac{d^\alpha f(x)}{dx^\alpha} &= \lim_{h\to 0} L\left(\frac{\tan(x+h)-\tan x}{(x+h)-x}\right)\\ &= \lim_{h\to 0} L\left(\frac{\frac{\tan x+\tanh}{1-\tan x\tanh}-\tan x}{h}\right)\\ &= \lim_{h\to 0}\left(\frac{\frac{\sec^2 h(1-\tan x\tanh)-(\tan x+\tanh)(-\tan x\sec^2 h}{(1-\tan x\tanh)^2}}{1}\right)\\ &= 1+\tan^2 x\\ &= \sec^2 x\end{aligned}$$

$$\frac{d^{\alpha} f(x)}{dx^{\alpha}} = \begin{matrix} \lim \\ h \to 0 \end{matrix} L\left( \frac{\cot(x+h) - \cot x}{(x+h) - x} \right)$$

$$= \begin{matrix} \lim \\ h \to 0 \end{matrix} L\left( \frac{\dfrac{\cos x \cosh - \sin xxinh}{\sin x \cosh + \cos x \sinh} - \cot x}{h} \right)$$

$$= \begin{matrix} \lim \\ h \to 0 \end{matrix} \left( \frac{\dfrac{-\sin^2 x - \cos^2 x}{\sin^2 x}}{1} \right)$$

$$= -(1 + \cot^2 x)$$

$$= \cos ec^2 x$$

### 3.4 Equaltiy for Logarithmic Functions

The derivative of logarithmic function is

$$\frac{d \log_a f(x)}{dx} = \frac{\dfrac{df(x)}{dx}}{f(x)} \ln a$$

For example, f(x)=lnx and the derivative is

$$\frac{d \ln x}{dx} = \frac{1}{x}$$

The FOD for α=1 is

$$f^{\alpha}(x) = \begin{matrix} \lim \\ h \to 0 \end{matrix} L\left( \frac{\ln(x+h) - \ln x}{(x+h) - x} \right)$$

$$= \frac{\dfrac{1}{x+h}}{1} = \frac{1}{x}$$

## 4 FOD for Real Order

In the case of real order of FOD, the definition can be revised as in Definition 3.

***Definion 3.*** *Assume that f(x):R→R is a function α∈R and L(.) is a L'Hospital process. FOD can be as follows:*

$$f^{(\alpha)}(x) = \begin{matrix} \lim \\ h \to 0 \end{matrix} L\left( \frac{f^{\alpha}(x+h) - f^{\alpha}(x)}{(x+h)^{\alpha} - x^{\alpha}} \right)$$

$$= \begin{matrix} \lim \\ h \to 0 \end{matrix} \frac{f'(x+h) f(x+h)^{\alpha-1}}{(x+h)^{\alpha-1}}$$

This definition can be applied to main function groups such polinomial, exponential, trigonometric and logarithmic. The obtained results can be expressed theorems.

***Theorem 3:****Assume a polinomial function is f(x)=x*$^n$*,* $\alpha = \dfrac{\beta}{\delta}$ *and δ≠0. Then the FOD of f(x) can be expressed as follows.*

$$f^{(\alpha)}(x) = \frac{nx^{n-1}\sqrt[\delta]{x^{n(\beta-\delta)}}}{\sqrt[\delta]{x^{\beta-\delta}}}$$

**Proof:** Assume that f(x) is a polinomial function f(x)=x$^n$, $\alpha = \dfrac{\beta}{\delta}$ and δ≠0. While the polinomial function is replaced in Definition 3, the following result can be obtained.

$$f^{(\alpha)}(x) = \lim_{h \to 0} \frac{f'(x+h)f(x+h)^{\alpha-1}}{(x+h)^{\alpha-1}}$$

It is known that the bionomial equation is

$$(x+h) = \binom{n}{0}x^n + \binom{n}{1}x^{n-1}h + \ldots + \binom{n}{n}h^n$$

and f(x)=x$^n$ can be replaced in definition

$$f^{(\alpha)}(x) = \lim_{h \to 0} \frac{\dfrac{d(x+h)^n}{dh}((x+h)^n)^{\alpha-1}}{(x+h)^{\alpha-1}}$$

$$= \lim_{h \to 0} \frac{\dfrac{d\left(x^n + nx^{n-1}h + \ldots + h^n\right)}{dh}(x+h)^{n\frac{\beta}{\delta}-1}}{(x+h)^{\frac{\beta}{\delta}-1}}$$

$$= \frac{nx^{n-1}\sqrt[\delta]{x^{n(\beta-\delta)}}}{\sqrt[\delta]{x^{\beta-\delta}}}$$

The equation in theorem was verified ■

***Theorem 4:*** *Assume an exponential function f(x)=x*$^n$*,* $\alpha = \dfrac{\beta}{\delta}$ *and δ≠0. Then the FOD of f(x) can be expressed as follows.*

$$f^{(\alpha)}(x) = \frac{e^x\sqrt[\delta]{e^{(\beta-\delta)x}}}{\sqrt[\delta]{x^{\beta-\delta}}}$$

**Proof:** Assume that an exponential function is f(x)=e$^x$, order is $\alpha = \dfrac{\beta}{\delta}$ and δ≠0. When an exponential function is substituted in Definition 3, the following FOD result can be obtained.

$$f^{(\alpha)}(x) = \lim_{h \to 0} \frac{f'(x+h) f(x+h)^{\alpha-1}}{(x+h)^{\alpha-1}}$$

$f(x)=e^x$ and its FOD is

$$\begin{aligned} f^{(\alpha)}(x) &= \lim_{h \to 0} \frac{\frac{de^{x+h}}{dh} (e^{x+h})^{\alpha-1}}{(x+h)^{\alpha-1}} \\ &= \lim_{h \to 0} \frac{\frac{e^x e^h}{1} \left(e^{x+h}\right)^{\frac{\beta}{\delta}-1}}{(x+h)^{\frac{\beta}{\delta}-1}} \\ &= \frac{e^x \sqrt[\delta]{e^{\beta-\delta}}}{\sqrt[\delta]{x^{\beta-\delta}}} \end{aligned}$$

∎

***Theorem 5**: Assume that f(x) is a trigonometric function, fractional order is* $\alpha = \frac{\beta}{\delta}$ *and* $\delta \neq 0$*. The FOD of f(x)=sinx is as follows.*

$$f^{(\alpha)}(x) = \frac{\cos x \sqrt[\delta]{(\sin x)^{\beta-\delta}}}{\sqrt[\delta]{x^{\beta-\delta}}}$$

**Proof:** Assume that a trigonometric function is f(x)=sinx, fractional order is $\alpha = \frac{\beta}{\delta}$ and $\delta \neq 0$. The FOD of sinx is follows.

$$\begin{aligned} f^{(\alpha)}(x) &= \lim_{h \to 0} \frac{\frac{d \sin(x+h)}{dh} (\sin(x+h))^{\alpha-1}}{(x+h)^{\alpha-1}} \\ &= \lim_{h \to 0} \frac{\frac{\cos(x+h)}{1} \left(\sin(x+h)\right)^{\frac{\beta}{\delta}-1}}{(x+h)^{\frac{\beta}{\delta}-1}} \\ &= \frac{\cos x \sqrt[\delta]{(\sin x)^{\beta-\delta}}}{\sqrt[\delta]{x^{\beta-\delta}}} \end{aligned}$$

∎

The FODs of the remaining trigonometric functions are seen in Table 3.

**Table 3.** FODs of trigonometric functions.

| Function | Fractional Order Derivative ($\alpha = \frac{\beta}{\delta}$ ve δ≠0) |
|---|---|
| sinx | $f^{(\alpha)} = \cos x \sqrt[\delta]{\left(\frac{\sin x}{x}\right)^{\beta-\delta}}$ |
| cosx | $f^{(\alpha)} = -\sin x \sqrt[\delta]{\left(\frac{\cos x}{x}\right)^{\beta-\delta}}$ |
| tanx | $f^{(\alpha)} = (1+\tan^2 x) \sqrt[\delta]{\left(\frac{\tan x}{x}\right)^{\beta-\delta}}$ |
| cotx | $f^{(\alpha)} = -(1+\cot^2 x) \sqrt[\delta]{\left(\frac{\cot x}{x}\right)^{\beta-\delta}}$ |
| secx | $f^{(\alpha)} = \sec x \tan x \sqrt[\delta]{\left(\frac{\sec x}{x}\right)^{\beta-\delta}}$ |
| cosecx | $f^{(\alpha)} = -\cos ecx \cot x \sqrt[\delta]{\left(\frac{\cos ecx}{x}\right)^{\beta-\delta}}$ |

***Theorem 6****: Assume thata logarithmic function f(x)=lnx, fractional order,* $\alpha = \frac{\beta}{\delta}$ *order is δ≠0. The following equation for FOD is sound and complete.*

$$f^{(\alpha)}(x) = \frac{1}{x}\sqrt[\delta]{\frac{(\ln x)^{\beta-\delta}}{x^{\beta-\delta}}}$$

**Proof:** Assume that a logarithmic function is f(x)=lnx, fractional order $\alpha = \frac{\beta}{\delta}$ and δ≠0. The fractional order derivative of logarithmic function can be expressed as in the following equation.

$$f^{(\alpha)}(x) = \lim_{h \to 0} \frac{\frac{d\ln(x+h)}{dh}(\ln(x+h))^{\alpha-1}}{(x+h)^{\alpha-1}}$$

$$= \lim_{h \to 0} \frac{\frac{1}{(x+h)}\left(\ln(x+h)\right)^{\frac{\beta}{\delta}-1}}{(x+h)^{\frac{\beta}{\delta}-1}}$$

$$= \frac{\sqrt[\delta]{(\ln x)^{\beta-\delta}}}{x\sqrt[\delta]{x^{\beta-\delta}}} = \frac{1}{x}\sqrt[\delta]{\frac{(\ln x)^{\beta-\delta}}{x^{\beta-\delta}}}$$

∎

## 5 Properties of FODs

The meaning of derivative is the rate of change in the dependent variable versus the changes in the independent variables. The applications of FODs can be handled by using monotonic increasing and decreasing functions.

***Theorem 7.*** *Assume that f(x) is a positive monotonic increasing function, $\alpha,\lambda \in R^+$ and $\alpha \le \lambda$. In this case, $f^{(\alpha)}(x) \le f^{(\lambda)}(x)$.*

**Proof.** f(x) is a monotonic increasing function, so $f(x_i) \le f(x_j)$ for $x_i \le x_j$ $\forall x_i, x_j \in R$. Then $f^\alpha(x_i+h) \le f^\lambda(x_j+h)$ and $f^\alpha(x_i) \le f^\lambda(x_j)$. This case implies that $f^\alpha(x_i+h)-f^\alpha(x_i) \le f^\lambda(x_j+h)-f^\lambda(x_j)$ and $(x_i+h)^\alpha - x_i^\alpha \le (x_j+h)^\lambda - x_j^\lambda$. So, this implies the following quotient.

$$\lim_{h\to 0} L\left(\frac{f^\alpha(x_i+h)-f^\alpha(x_i)}{(x_i+h)^\alpha - x_i^\alpha}\right) \le \lim_{h\to 0} L\left(\frac{f^\lambda(x_j+h)-f^\lambda(x_j)}{(x_j+h)^\lambda - x_j^\lambda}\right) \blacksquare$$

***Theorem 8.*** *Assume that f(x) is a positive monotonic decreasing function, $\alpha,\lambda \in R^+$ and $\alpha \le \lambda$. In this case, $f^{(\alpha)}(x) \ge f^{(\lambda)}(x)$.*

**Proof.** f(x) is a monotonic increasing function, so $f(x_i) \ge f(x_j)$ for $x_i \le x_j$ $\forall x_i, x_j \in R$. $f^\alpha(x_i+h) \ge f^\lambda(x_j+h)$ and $f^\alpha(x_i) \ge f^\lambda(x_j)$. This case implies that $f^\alpha(x_i+h)-f^\alpha(x_i) \ge f^\lambda(x_j+h)-f^\lambda(x_j)$ and $(x_i+h)^\alpha - x_i^\alpha \le (x_j+h)^\lambda - x_j^\lambda$. The term $f^\alpha(x_i+h)-f^\alpha(x_i)$ is divided by a smaller term such as $(x_i+h)^\alpha - x_i^\alpha$ than the term $(x_j+h)^\lambda - x_j^\lambda$ divides $f^\lambda(x_j+h)-f^\lambda(x_j)$. So, this implies the following quotient.

$$\lim_{h\to 0} L\left(\frac{f^\alpha(x_i+h)-f^\alpha(x_i)}{(x_i+h)^\alpha - x_i^\alpha}\right) \ge \lim_{h\to 0} L\left(\frac{f^\lambda(x_j+h)-f^\lambda(x_j)}{(x_j+h)^\lambda - x_j^\lambda}\right) \blacksquare$$

***Theorem 9.*** *Assume that f(x) is a positive monotonic increasing function, $\alpha,\lambda \in R^-$ and $\alpha \le \lambda$. In this case, $f^{(\alpha)}(x) \le f^{(\lambda)}(x)$.*

**Proof.** f(x) is a positive monotonic increasing function, so $f(x_i) \le f(x_j)$ for $x_i \le x_j$ $\forall x_i, x_j \in R$. Then $f^\alpha(x_i+h) \le f^\lambda(x_j+h)$ and $f^\alpha(x_i) \le f^\lambda(x_j)$, $|\alpha| \ge |\lambda|$.

$$A = \lim_{h\to 0} L\left(\frac{f^\alpha(x_i+h)-f^\alpha(x_i)}{(x_i+h)^\alpha - x_i^\alpha}\right)$$

$$= \frac{\dfrac{d}{dh}\dfrac{\left(f^{|\alpha|}(x)-f^{|\alpha|}(x+h)\right)}{f^{|\alpha|}(x)\,f^{|\alpha|}(x+h)}}{\dfrac{d}{dh}\dfrac{x^{|\alpha|}-(x+h)^{|\alpha|}}{x^{|\alpha|}(x+h)^{|\alpha|}}}$$

$$= \frac{f'(x)(x)^{|\alpha|+1}}{f^{|\alpha|+1}(x)}$$

and

$$B = \lim_{h \to 0} L\left( \frac{f^{\lambda}(x_i + h) - f^{\lambda}(x_i)}{(x_i + h)^{\lambda} - x_i^{\lambda}} \right)$$

$$= \frac{\frac{d}{dh} \frac{\left(f^{|\lambda|}(x) - f^{|\lambda|}(x+h)\right)}{f^{|\lambda|}(x) f^{|\lambda|}(x+h)}}{\frac{d}{dh} \frac{x^{|\lambda|} - (x+h)^{|\lambda|}}{x^{|\lambda|}(x+h)^{|\lambda|}}} .$$

$$= \frac{f'(x)(x)^{|\lambda|+1}}{f^{|\lambda|+1}(x)}$$

$f^{|\lambda|+1}(x) \leq f^{|\alpha|+1}(x)$ and $x^{|\lambda|+1} \leq x^{|\alpha|+1}$. This case implies the following inequality.

$$A \leq B \Rightarrow \frac{f'(x)(x)^{|\alpha|+1}}{f^{|\alpha|+1}(x)} \leq \frac{f'(x)(x)^{|\lambda|+1}}{f^{|\lambda|+1}(x)}$$

$$\Rightarrow \frac{(x)^{|\alpha|+1}}{f^{|\alpha|+1}(x)} \leq \frac{(x)^{|\lambda|+1}}{f^{|\lambda|+1}(x)} .$$

This is the completeness of proof ■

***Theorem 10.*** *Assume that f(x) is a positive monotonic decreasing function,* $\alpha, \lambda \in R^-$ *and* $\alpha \leq \lambda$*. In this case,* $f^{(\alpha)}(x) \geq f^{(\lambda)}(x)$*.*

**Proof.** g(x) is a positive monotonic decreasing function, so $g(x_i) \geq g(x_j)$ for $x_i \leq x_j$ $\forall x_i, x_j \in R$. Then $g^{\alpha}(x_i+h) \geq g^{\lambda}(x_j+h)$ and $g^{\alpha}(x_i) \leq g^{\lambda}(x_j)$, $|\alpha| \geq |\lambda|$. Assume that f(x) is a positive monotonic increasing function, then $g(x) = \frac{1}{f(x)}$ is a positive monotonic decreasing function.

$$A = \lim_{h \to 0} L\left( \frac{g^{\alpha}(x_i + h) - g^{\alpha}(x_i)}{(x_i + h)^{\alpha} - x_i^{\alpha}} \right)$$

$$= \frac{\frac{d}{dh}\left(f^{|\alpha|}(x+h) - f^{|\alpha|}(x)\right)}{\frac{d}{dh} \frac{x^{|\alpha|} - (x+h)^{|\alpha|}}{x^{|\alpha|}(x+h)^{|\alpha|}}}$$

$$= -\frac{f'(x) f^{|\alpha|-1}(x)}{x^{3|\alpha|-1}}$$

and

$$
\begin{aligned}
B &= \lim_{h \to 0} L\left( \frac{g^{\lambda}(x_i + h) - g^{\lambda}(x_i)}{(x_i + h)^{\lambda} - x_i^{\lambda}} \right) \\
&= \frac{\frac{d}{dh}\left( f^{|\lambda|}(x+h) - f^{|\lambda|}(x) \right)}{\frac{d}{dh}\frac{x^{|\lambda|} - (x+h)^{|\lambda|}}{x^{|\lambda|}(x+h)^{|\lambda|}}} \\
&= -\frac{f'(x) f^{|\lambda|-1}(x)}{x^{3|\lambda|-1}}
\end{aligned}
.
$$

$f^{|\lambda|-1}(x) \leq f^{|\alpha|-1}(x)$ and $x^{3|\lambda|-1} \leq x^{3|\alpha|-1}$. This case implies the following inequality.

$$
\begin{aligned}
A \leq B &\Rightarrow -\frac{f'(x) f^{|\alpha|-1}(x)}{x^{3|\alpha|-1}} \geq -\frac{f'(x) f^{|\lambda|-1}(x)}{x^{3|\lambda|-1}} \\
&\Rightarrow -\frac{f'(x) f^{|\alpha|-1}(x)}{x^{3|\alpha|-1}} \geq -\frac{f'(x) f^{|\lambda|-1}(x)}{x^{3|\lambda|-1}}
\end{aligned}
.
$$

This is the completeness of proof ■

The fractional order derivative can be characterized for any functions.

***Theorem 11.*** *Assume that f(x) is a function such as f:R→R and $\alpha \in R$. then $f^{(\alpha)}(x)$ is a function of complex variables.*

**Proof.** Assume that $\alpha = \frac{\beta}{\delta}$ and $\delta \neq 0$. If $f(x) \geq 0$, the results obtained in Theorem 1, Theorem 2, Theorem 3 and Theorem 4 are the cases. The fractional order derivative of f(x) is

$$
\begin{aligned}
f^{(\alpha)} &= \frac{f'(x) f^{\alpha-1}(x)}{x^{\alpha-1}} \\
&= f'(x) \sqrt[\delta]{\frac{f^{\beta-\delta}(x)}{x^{\beta-\delta}}} \\
&= f'(x) \sqrt[\delta]{\left( \frac{f(x)}{x} \right)^{\beta-\delta}}
\end{aligned}
$$

If the FOD is a function of complex variables, then $f^{(\alpha)}(x) = g(x) + ih(x)$ where i is the $i = \sqrt{-1}$.

If $f(x) < 0$, there will be two cases:

**Case 1:** Assume that $\delta$ is odd.

If $\left( \frac{f(x)}{x} \right)^{\beta-\delta} \geq 0$ or $\left( \frac{f(x)}{x} \right)^{\beta-\delta} < 0$,

then the obtained function $f^{(\alpha)}(x)$ is a real function and h(x)=0 for both cases. Since the multiplication of any negative number in odd steps yields a negative number.

**Case 2:** Assume that δ is even.

If $\left(\frac{f(x)}{x}\right)^{\beta-\delta} \geq 0$,

then h(x)=0 and $f^{(\alpha)}(x)$ is a real function.

If $\left(\frac{f(x)}{x}\right)^{\beta-\delta} < 0$,

then the multiplication of any number in even steps yields a positive number for real numbers. However, it yields a negative result for complex numbers, so, h(x)≠0. This means that $f^{(\alpha)}(x)$ is a complex function.

In fact, $f^{(\alpha)}(x)$ is a complex function for both cases. The h(x)=0 for some situations ∎

## 6 Implementations of New Approach for Fractional Order Derivative

The implementations were handled for positive monotonic increasing/decreasing functions. So there are four cases.

a) Assume that f(x) is a positive monotonic increasing function in [0,∞). $f(x)=x^2+3x+4$ is a function obeys this case. The FODs of f(x) for α=-0.5, α=0.5, α=1, α=1.5 and α=2.5 is seen in Fig.4. The satisfaction of claim in Theorem 7 is illustrated in Fig.4.

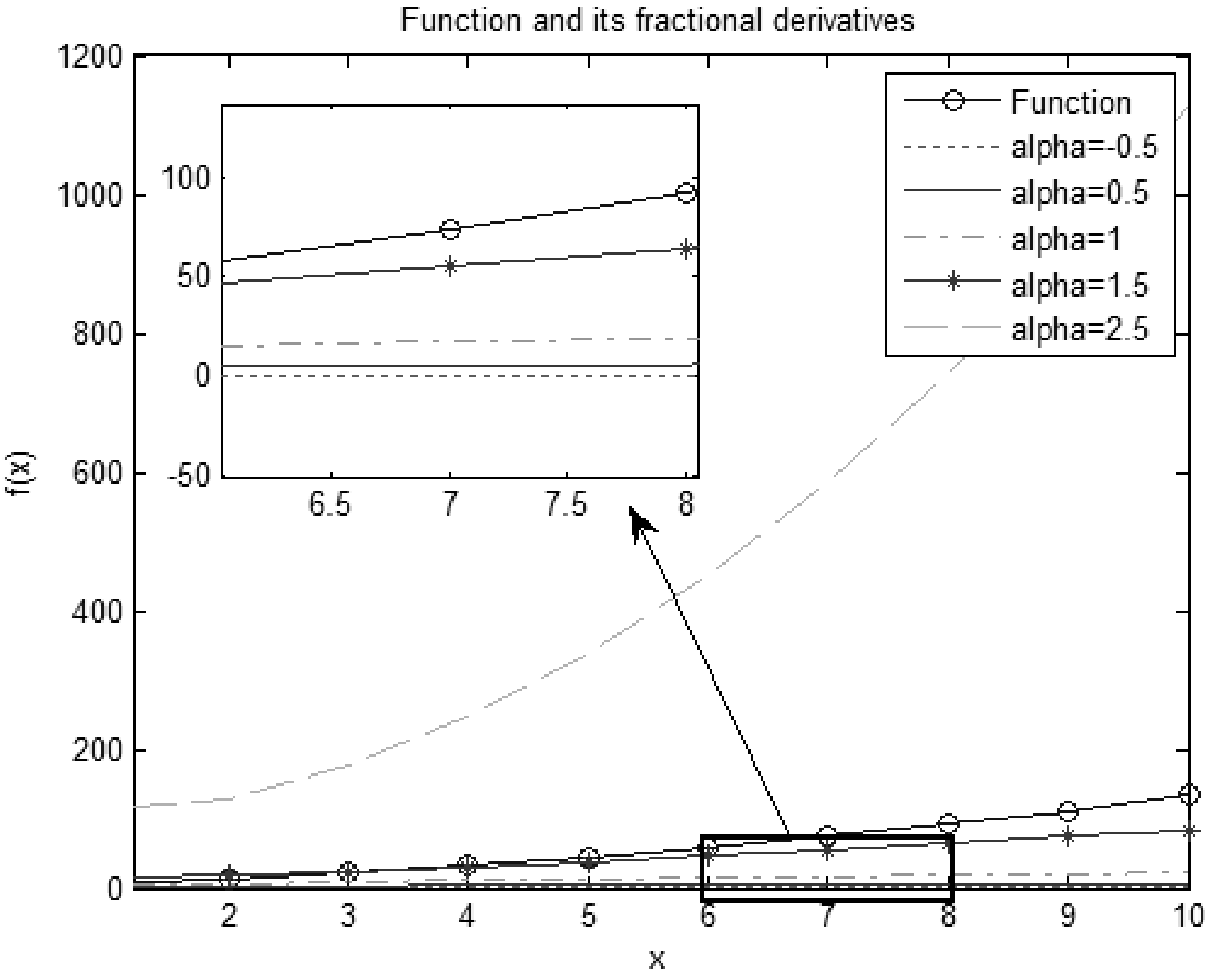

**Figure 4.** The FODs of $f(x)=x^2+3x+4$ and function itself.

b) Assume that f(x) is a positive monotonic decreasing function in [0,1]. $f(x)=\frac{1}{x^2+3x+4}$ is a function obeys this case. The FODs of f(x) for α=-0.5, α=0.5, α=1, α=1.5 and α=2.5 is seen in Fig.5. The satisfaction of claim in Theorem 8 is illustrated in Fig.5.

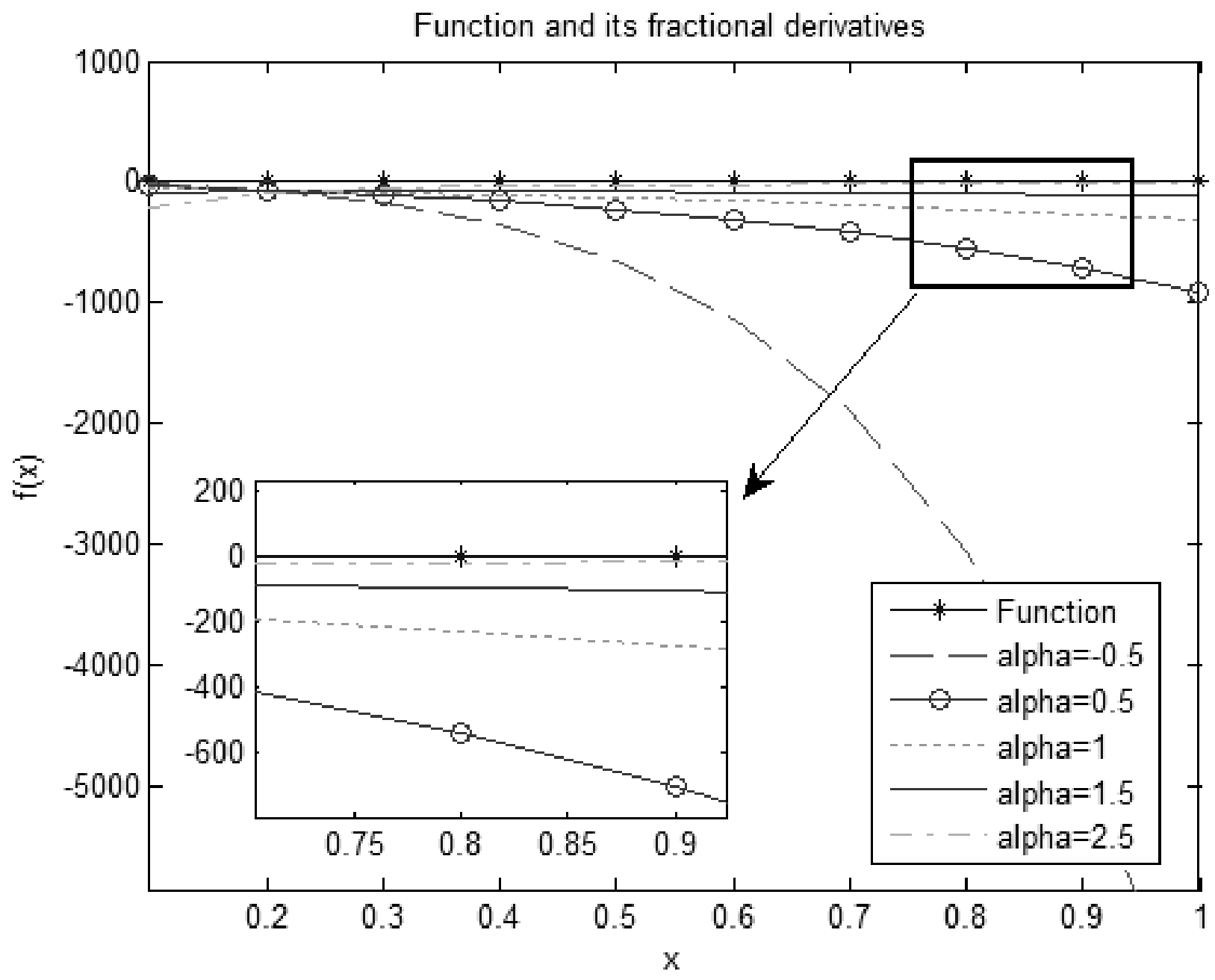

**Figure 5.** The FODs of $f(x) = \dfrac{1}{x^2 + 3x + 4}$ and function itself.

c) Assume that f(x) is a positive monotonic increasing function in [0.001,0.1]. $f(x)=x^2+3x+0.001$ is a function obeys this case. The FODs of f(x) for $\alpha$=-2.5, $\alpha$=-2, $\alpha$=-1.5, $\alpha$=-1 and $\alpha$=-0.5 is seen in Fig.6. The satisfaction of claim in Theorem 9 is illustrated in Fig.6.

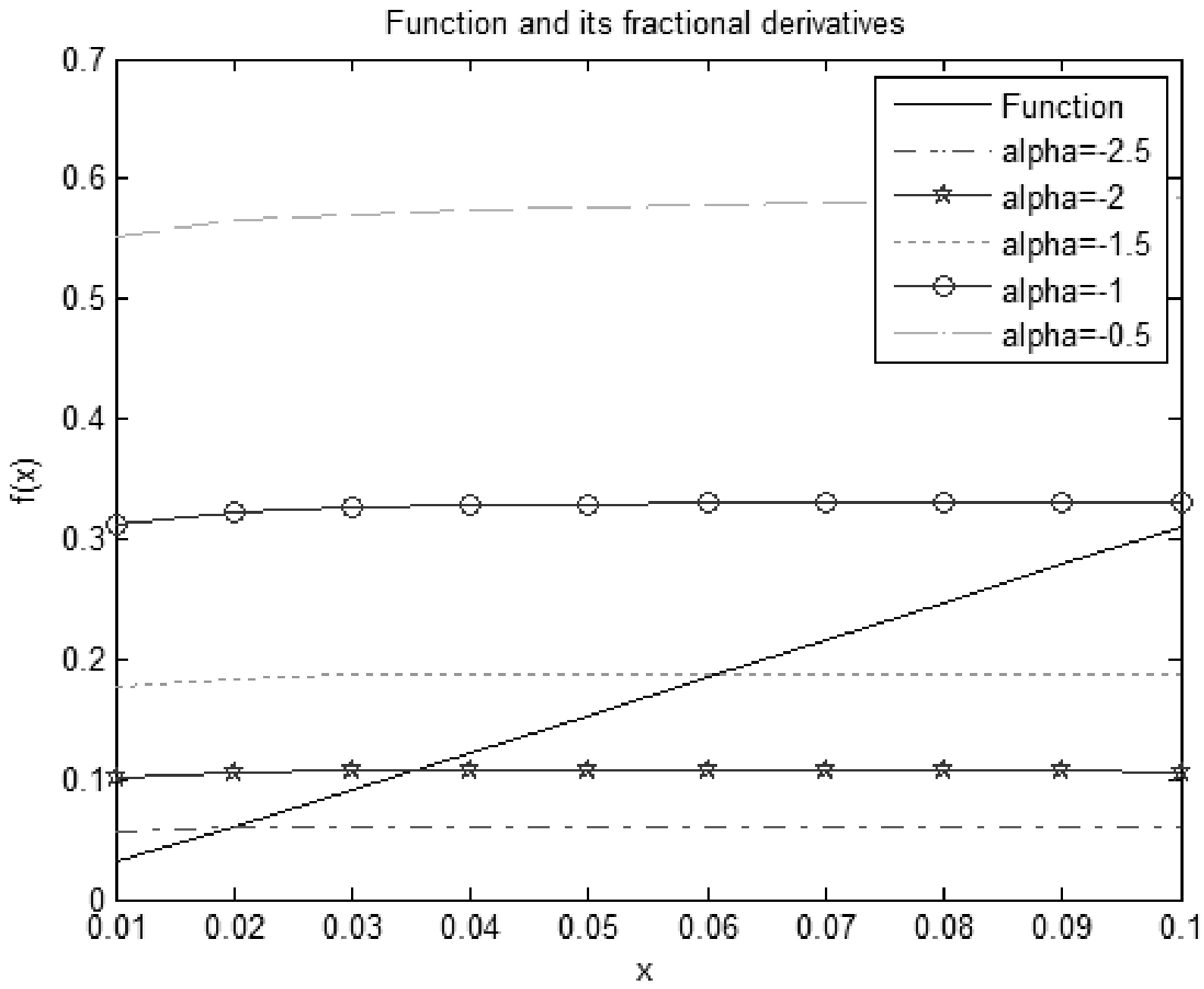

**Figure 6.** The fractional order derivatives of $f(x)=x^2+3x+0.001$ and function itself.

d) Assume that f(x) is a positive monotonic decreasing function in [0.001,0.1]. $f(x)=\dfrac{1}{x^2+3x+0.001}$ is a function obeys this case. The FODs of f(x) for α=-2.5, α=-2, α=-1.5, α=-1 and α=-0.5 is seen in Fig.7. The satisfaction of claim in Theorem 10 is illustrated in Fig.7.

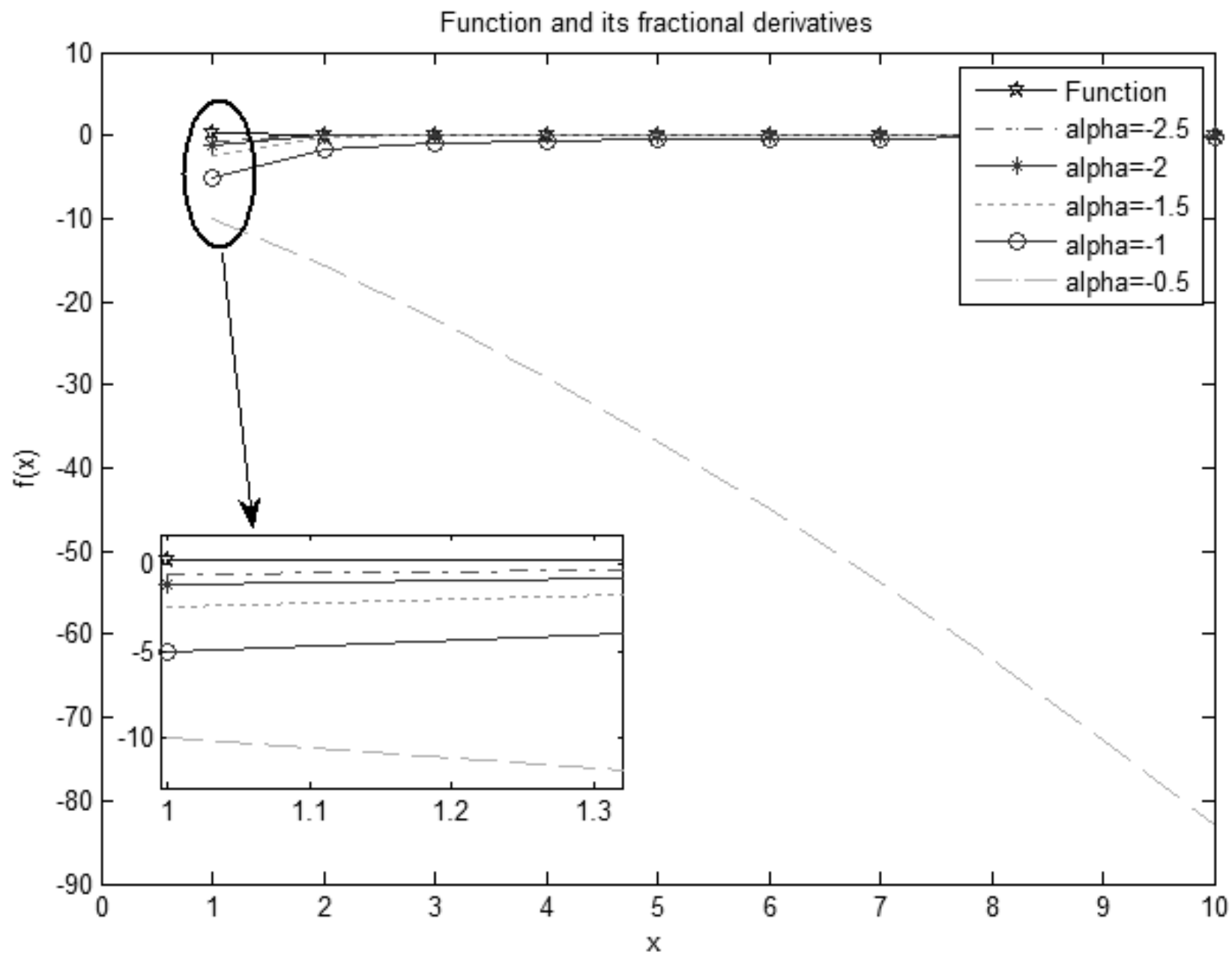

**Figure 7.** The FODs of $f(x)=\dfrac{1}{x^2+3x+0.001}$ and function itself.

The applications of FOD for logarithmic, trigonometric, exponential and polynomial functions were done in this section. The FOD of sin(x) is $f^{(\alpha)}(x)=\dfrac{\cos(x)(\sin(x))^{\alpha-1}}{x^{\alpha-1}}$. The graphics of FOD for α=0.5, α=1, α=1.5, α=2, α=2.5 and function itself are seen in Figure 8.

The FOD of cos(x) is $f^{(\alpha)}(x)=\dfrac{-\sin(x)(\cos(x))^{\alpha-1}}{x^{\alpha-1}}$. The graphics of FOD for α=0.5, α=1, α=1.5, α=2, α=2.5 and function itself are seen in Figure 9. The FOD of tan(x) is $f^{(\alpha)}(x)=\dfrac{(1+\tan^2(x))(\tan(x))^{\alpha-1}}{x^{\alpha-1}}$. The graphics of FOD for α=0.5, α=1, α=1.5, α=2, α=2.5 and function itself are seen in Figure 10.

The FOD of cot(x) is $f^{(\alpha)}(x)=\dfrac{-(1+\cot^2(x))(\cot(x))^{\alpha-1}}{x^{\alpha-1}}$. The graphics of FOD for α=0.5, α=1, α=1.5, α=2, α=2.5 and function itself are seen in Figure 11. The FOD of $f(x)=x^4-5x^3+x-2$ is $f^{(\alpha)}(x)=\dfrac{(4x^3-15x^2+1)(x^4-5x^3+x-2)^{\alpha-1}}{x^{\alpha-1}}$. The graphics of FOD for α=0.5, α=1, α=1.5, α=2, α=2.5 and function itself are seen in Figure 12.

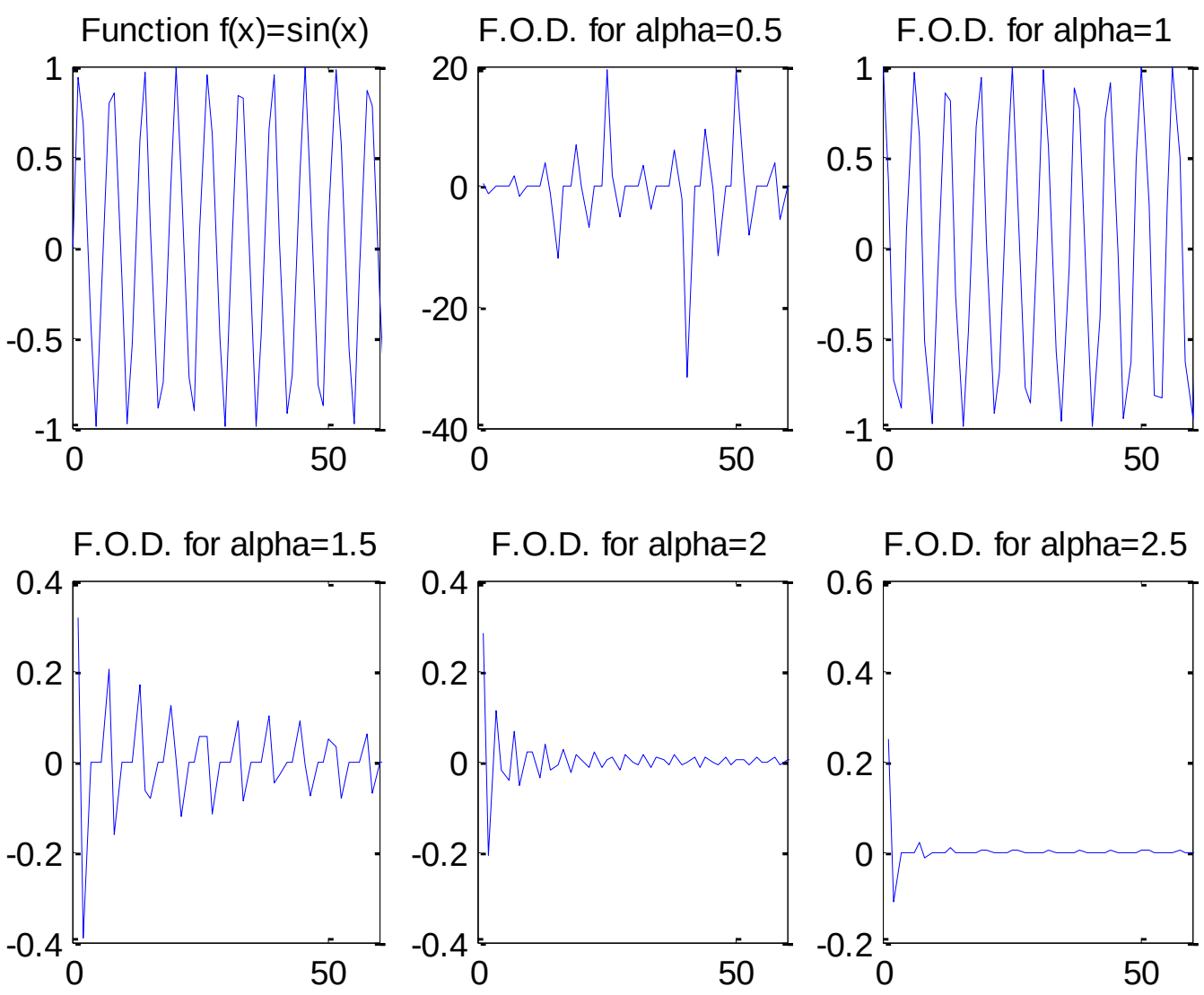

**Figure 8.** Graphics of FOD for α=0.5, α=1, α=1.5, α=2, α=2.5 and function itself.

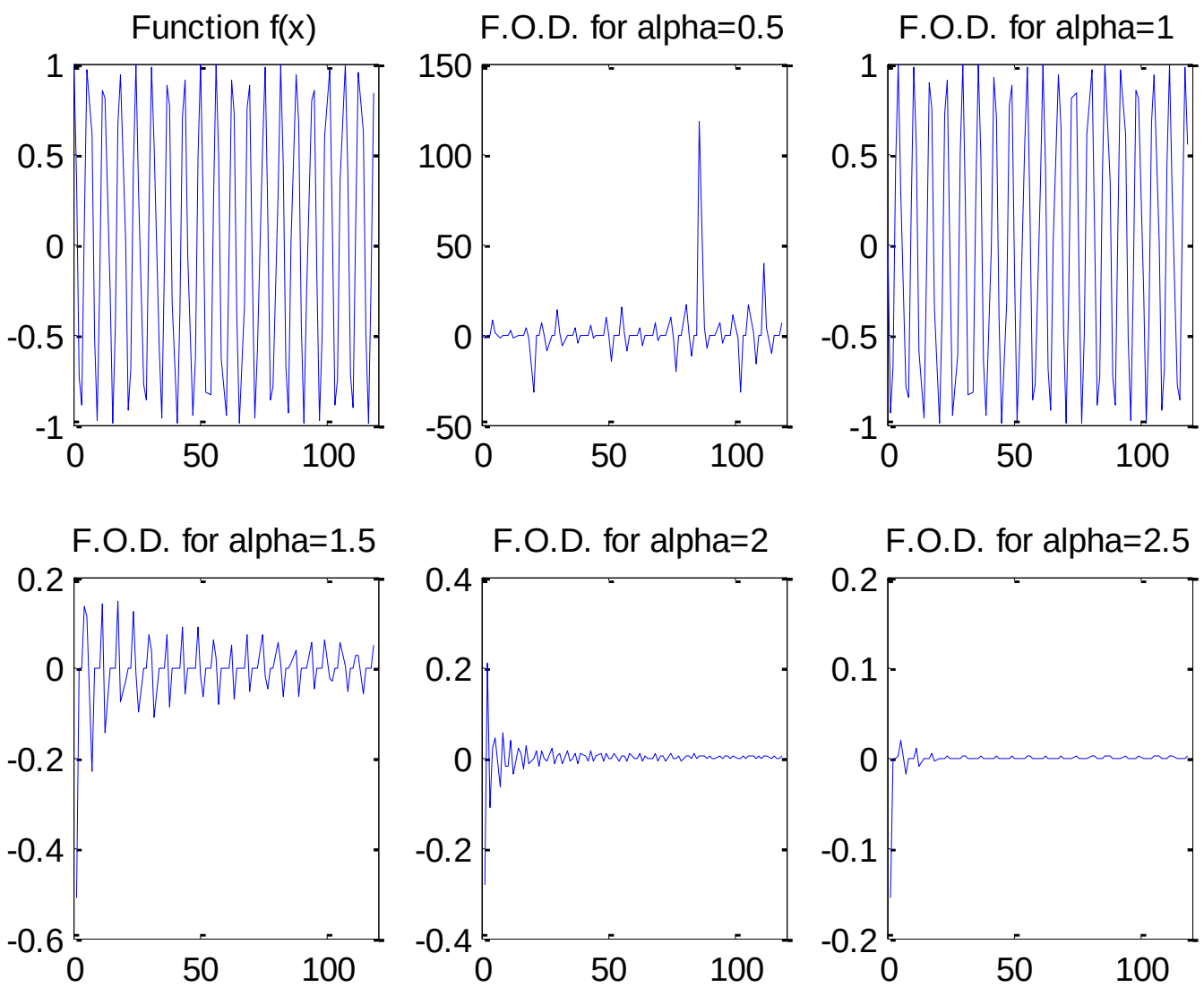

**Figure 9.** Graphics of FOD for α=0.5, α=1, α=1.5, α=2, α=2.5 and function itself (f(x)=cos(x)).

The FOD of $f(x)=\dfrac{1}{x^4-5x^3+x-2}$ is $f^{(\alpha)}(x)=\dfrac{\left(-\dfrac{4x^3-15x^2+1}{\left(x^4-5x^3+x-2\right)^2}\right)\left(\dfrac{1}{x^4-5x^3+x-2}\right)^{\alpha-1}}{x^{\alpha-1}}$.
The graphics of FOD for α=0.5, α=1, α=1.5, α=2, α=2.5 and function itself are seen in Figure 13.

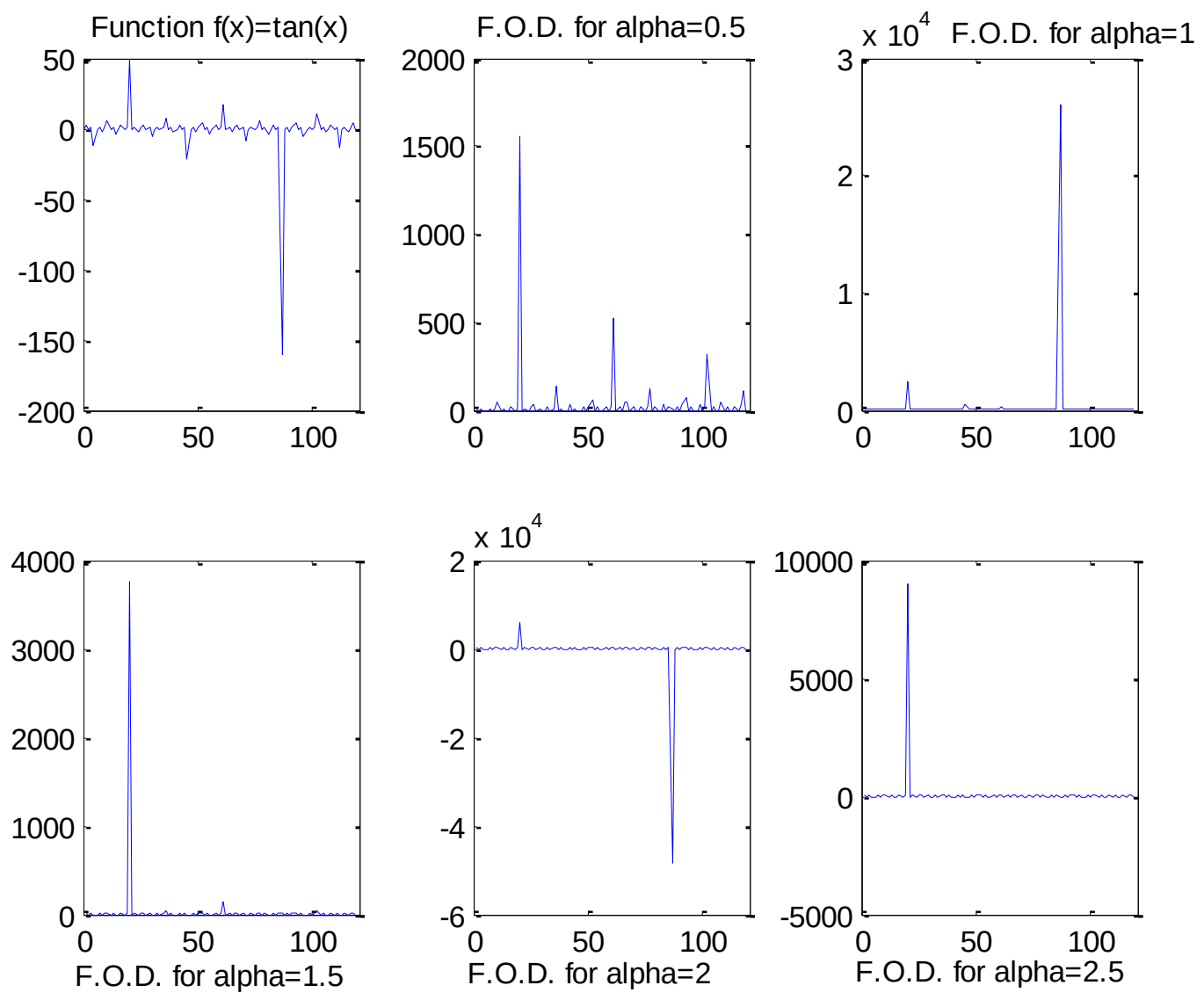

**Figure 10.** Graphics of FOD for α=0.5, α=1, α=1.5, α=2, α=2.5 and function itself (f(x)=tan(x)).

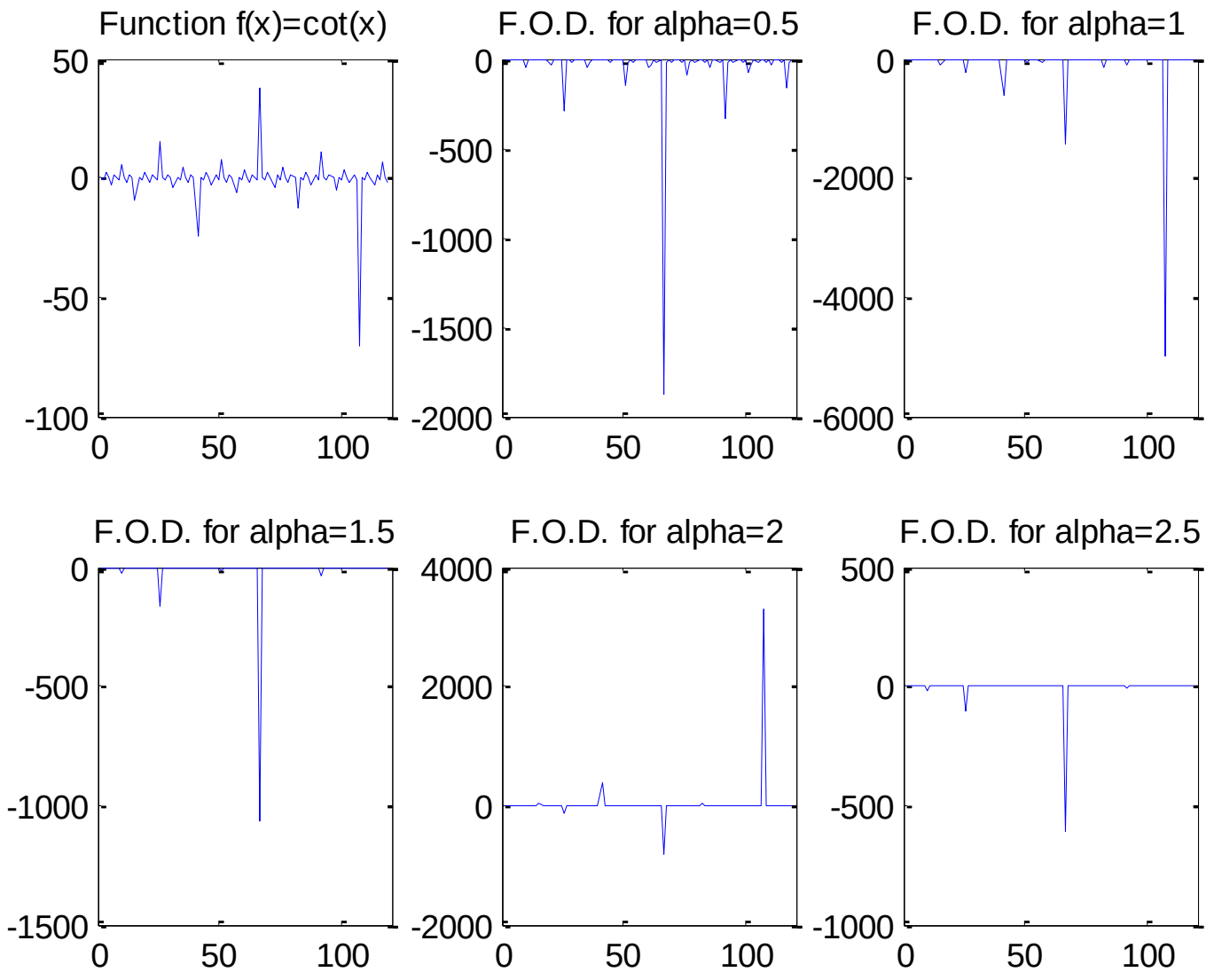

**Figure 11.** Graphics of FOD for α=0.5, α=1, α=1.5, α=2, α=2.5 and function itself (f(x)=cot(x)).

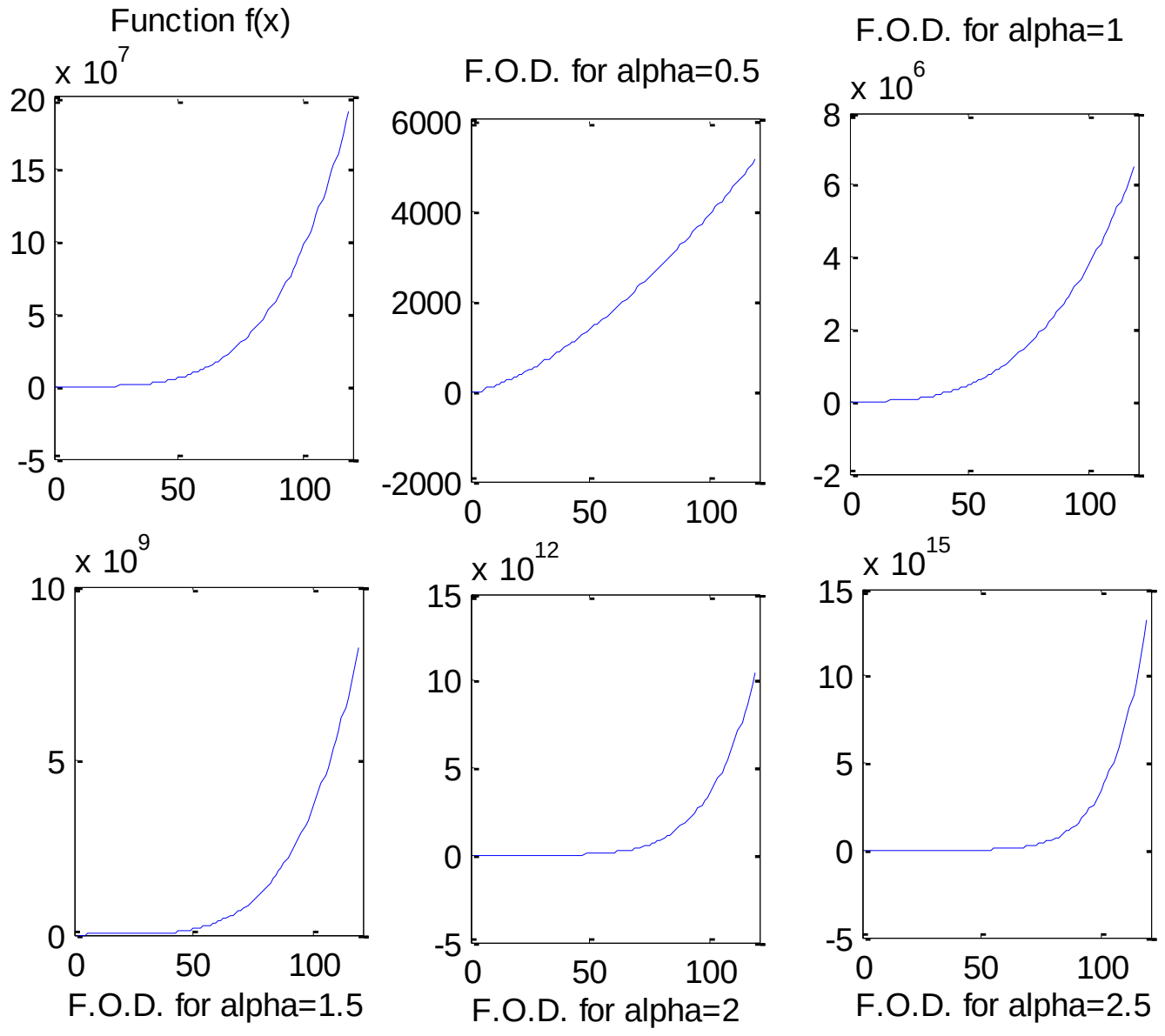

**Figure 12.** Graphics of FOD for α=0.5, α=1, α=1.5, α=2, α=2.5 and function itself (f(x)= $=x^4-5x^3+x-2$).

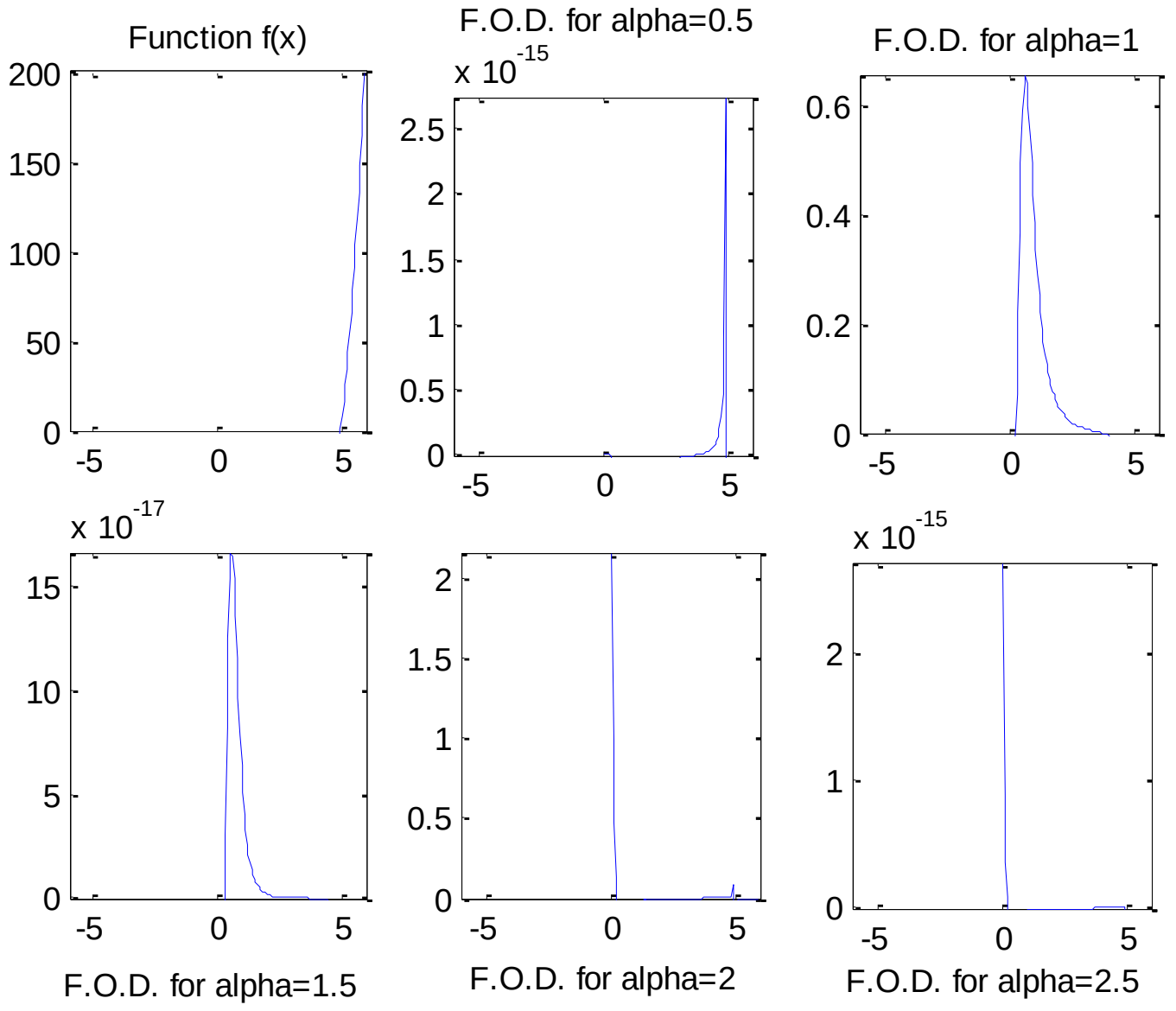

**Figure 13.** Graphics of FOD for α=0.5, α=1, α=1.5, α=2, α=2.5 and function itself ($f(x)=\dfrac{1}{x^4-5x^3+x-2}$).

The FOD of f(x)=ln(x) is $f^{(\alpha)}(x)=\dfrac{\left(\frac{1}{x}\right)(\ln(x))^{\alpha-1}}{x^{\alpha-1}}$. The graphics of FOD for α=0.5, α=1, α=1.5, α=2, α=2.5 and function itself are seen in Figure 14. The FOD of $f(x)=e^x$ is $f^{(\alpha)}(x)=\dfrac{e^x\left(e^x\right)^{\alpha-1}}{x^{\alpha-1}}$.

The graphics of FOD for α=0.5, α=1, α=1.5, α=2, α=2.5 and function itself are seen in Figure 15. The FOD of f(x)=$2^x$ is $f^{(\alpha)}(x) = \frac{2^x \ln(2)\left(2^x\right)^{\alpha-1}}{x^{\alpha-1}}$. The graphics of FOD for α=0.5, α=1, α=1.5, α=2, α=2.5 and function itself are seen in Figure 16.

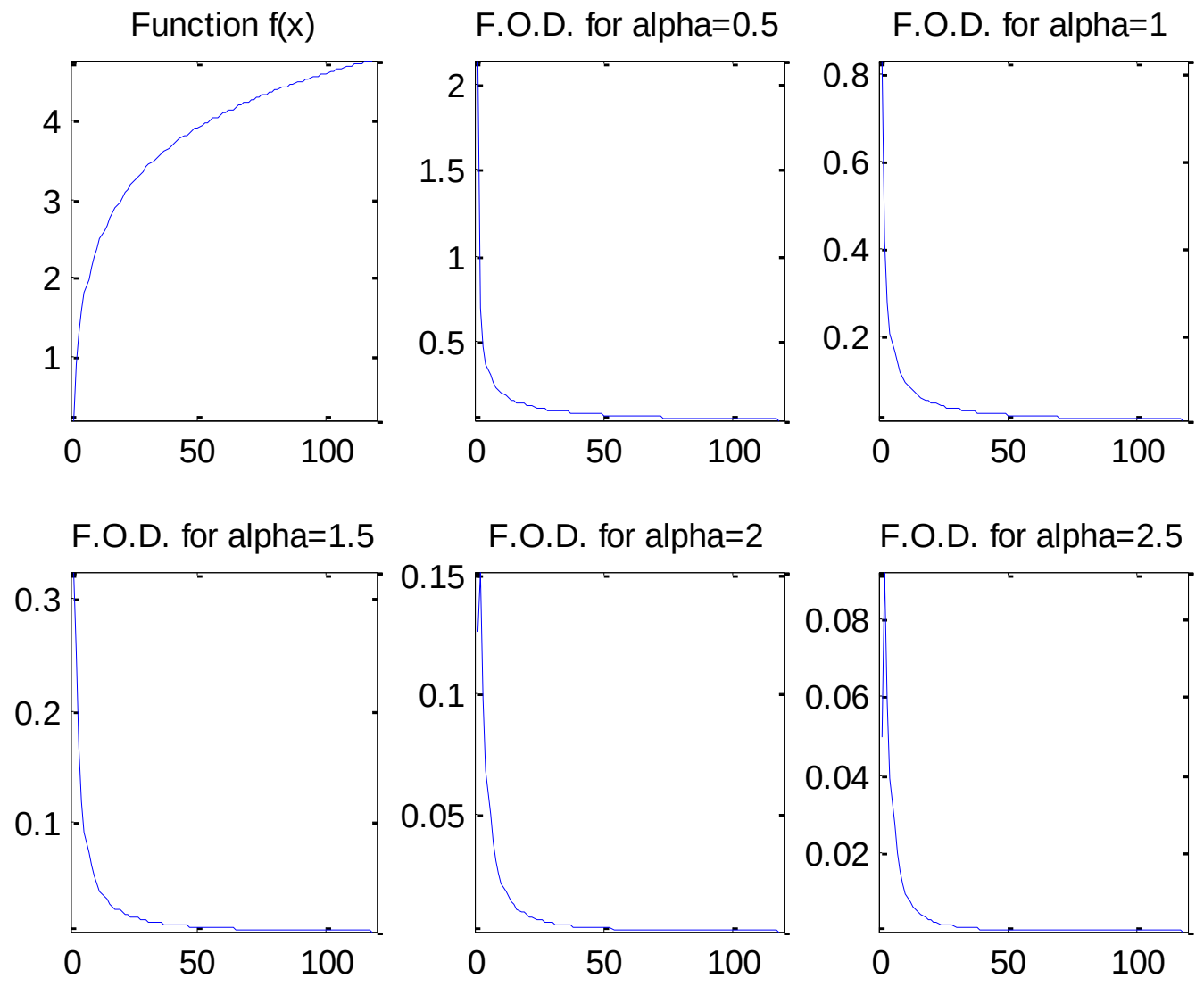

**Figure 14.** Graphics of FOD for α=0.5, α=1, α=1.5, α=2, α=2.5 and function itself (f(x)=ln(x)).

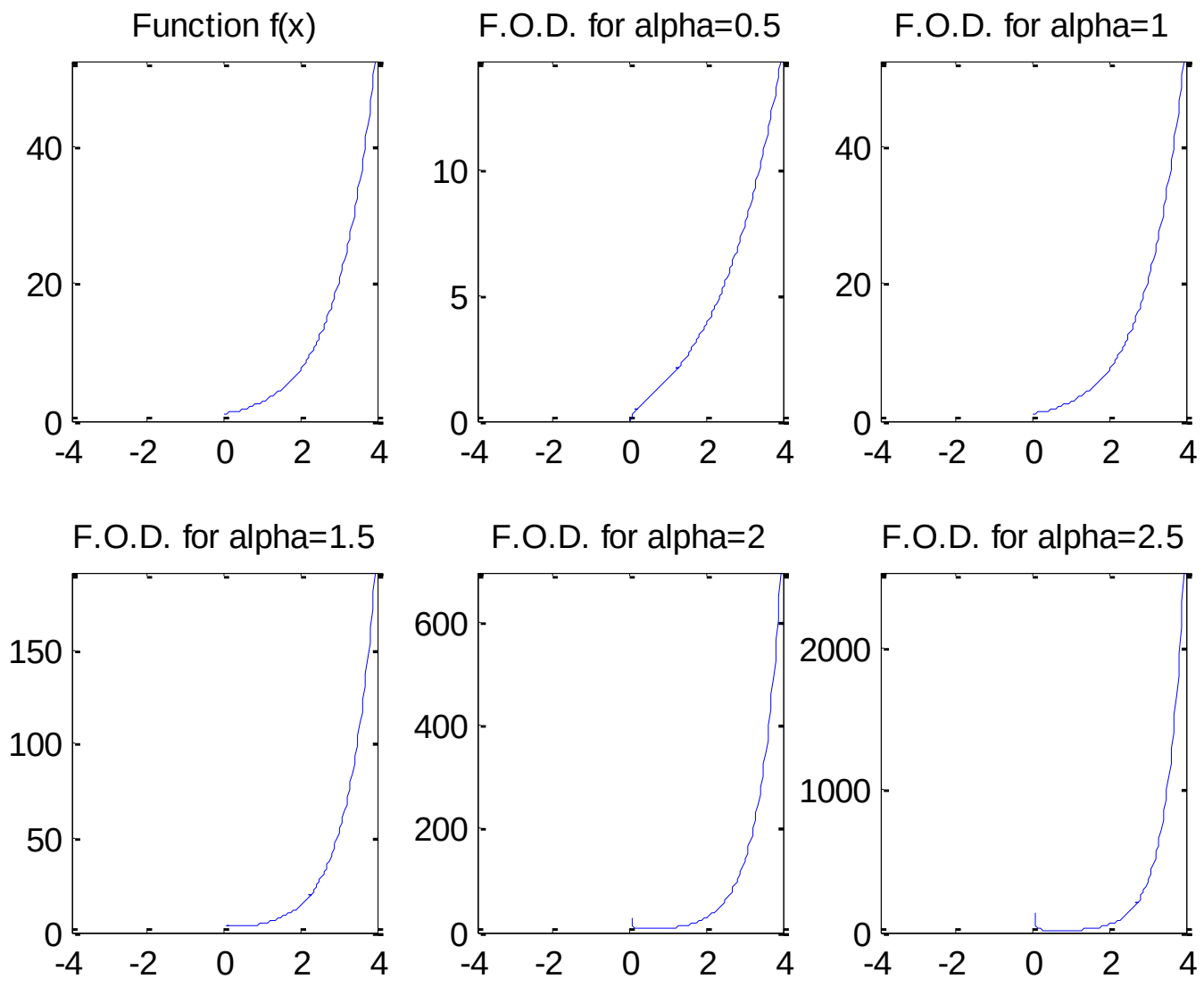

**Figure 15.** Graphics of FOD for α=0.5, α=1, α=1.5, α=2, α=2.5 and function itself ( $f(x) = e^x$ ).

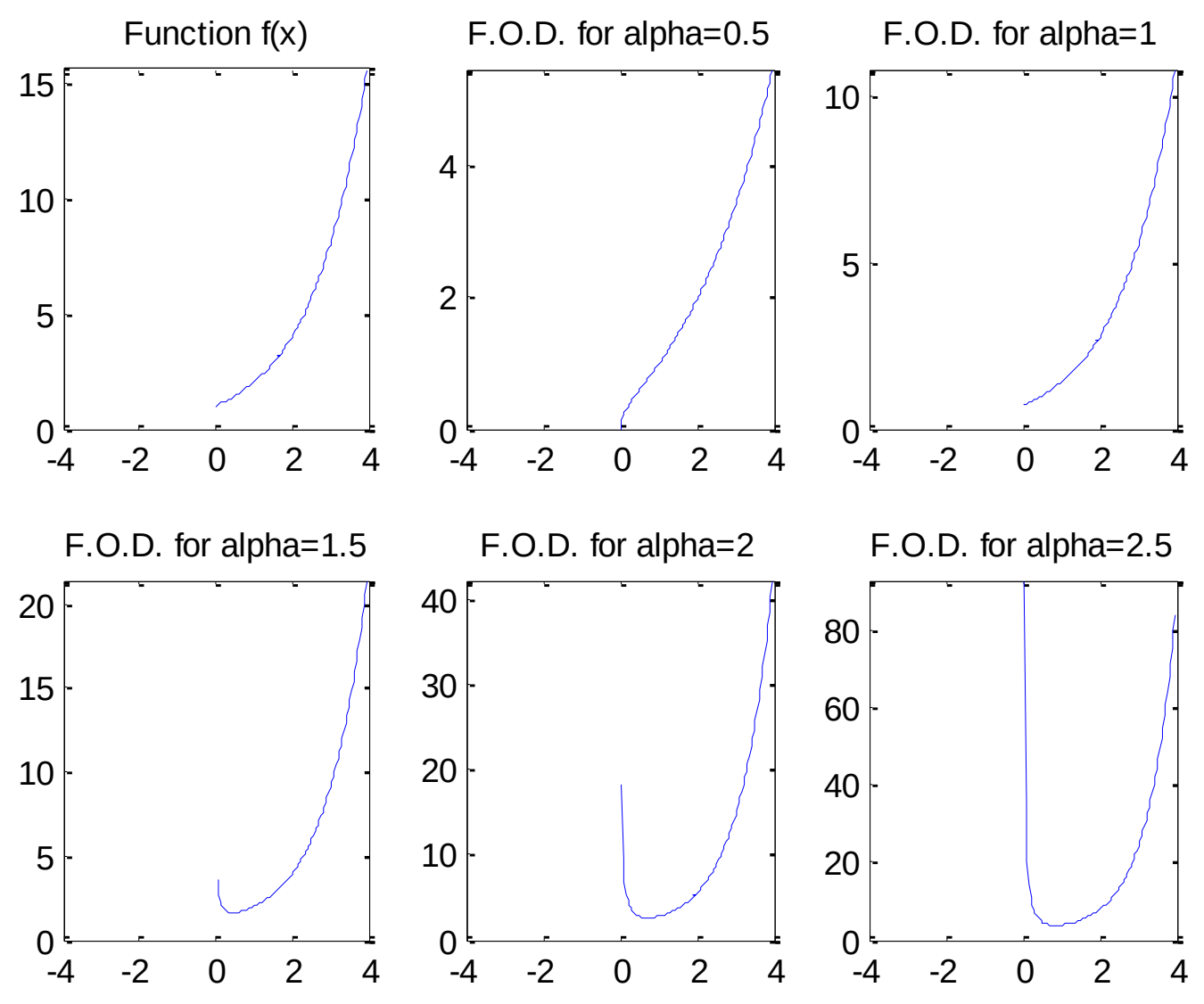

**Figure 16.** Graphics of FOD for α=0.5, α=1, α=1.5, α=2, α=2.5 and function itself ($f(x)=2^x$).

## 7 Conclusions

The old FOD methods have some important deficiencies, since they assumed the order of derivation is integer up to a specific derivation step. This caused involvement of deficiencies in the obtained formula. The derivation process can be redefined in the case of real order. The FODs for some important function groups (polynomial, exponential, trigonometric and logarithmic) were defined in this paper. While order is selected as 1, the obtained results are same as results of classical derivative. The properties of FODs handled in this study and finally the results of derivation are related to complex functions. However, this case has not been handled up to today. The results of derivation are a function of complex variables, or a complex number. This case is illustrated in this study.

The new definition for FOD causes some important conjectures:

1) The geometric meaning of FOD (the order is positive) can be analysed and which value of fractional order is applicable in practise. Both of these cases subject to future researches.
2) The geometric meaning of FOD (the order is negative) can be analysed and which value of fractional order is applicable in practise. Both of these cases subject to future researches.
3) The FOD can be analysed in case of complex order. The geometric and computational meaning must be searched in the future researches. For complex orders.
4) The integral of FOD is also a research area for the future studies. The answers of these questions cause the new research areas.

Ali KARCI
İnönü University
Faculty of Engineering
Department of Computer Engineering
44280, Malatya / Turkey
ali.karci@inonu.edu.tr; adresverme@gmail.com